\documentclass[a4paper]{paper}
\usepackage{hyperref}
\usepackage{amsmath}
\usepackage{amssymb,amsfonts,amsthm}
\usepackage{mathrsfs}
\usepackage{cleveref}
\usepackage{acro}
\usepackage[section]{placeins}
\usepackage{algorithm}
\usepackage[noend]{algpseudocode}
\usepackage{graphicx}
\usepackage{subcaption}
\usepackage[squaren,thickqspace]{SIunits}
\usepackage{booktabs}

\newcommand{\Real}{\ensuremath{\mathbb{R}}}

\newcommand{\Natural}{\ensuremath{\mathbb{N}}}
\newcommand{\domain}{\Omega}
\newcommand\datamanifold{\mathbb{M}}

\newcommand{\RecSpace}{\ensuremath{X}}
\newcommand{\DataSpace}{\ensuremath{Y}}
\newcommand{\ParamSpace}{\ensuremath{Z}}
\newcommand{\ProdSpace}{\ensuremath{U}}
\newcommand{\FeatureSpace}{\ensuremath{\mathbb{F}}}
\newcommand{\approxInv}[1]{#1^{\dagger}}
\newcommand{\learnedInv}[1]{#1^{\dagger}}
\DeclareMathOperator{\IdentityOp}{\ensuremath{\text{Id}}}
\DeclareMathOperator{\OpA}{\ensuremath{\mathcal{A}}}
\DeclareMathOperator{\OpB}{\ensuremath{\mathcal{B}}}
\DeclareMathOperator{\OpC}{\ensuremath{\mathcal{C}}}
\DeclareMathOperator{\OpF}{\ensuremath{\mathcal{F}}}
\DeclareMathOperator{\OpG}{\ensuremath{\mathcal{G}}}
\DeclareMathOperator{\OpK}{\ensuremath{\mathcal{K}}}
\DeclareMathOperator{\ForwardOp}{\ensuremath{\mathcal{T}}}
\DeclareMathOperator{\ForwardOpPseudoInv}{\ensuremath{\approxInv{\ForwardOp}}}
\DeclareMathOperator{\ForwardOpInvLearned}{\ensuremath{\learnedInv{\ForwardOp}_{\vparam}}}
\DeclareMathOperator{\RadonTransform}{\ensuremath{\mathcal{P}}}
\DeclareMathOperator{\LogLikelihood}{\ensuremath{\mathcal{L}}}
\DeclareMathOperator{\AffineOp}{\ensuremath{\mathcal{W}}}
\DeclareMathOperator{\NonLinOp}{\ensuremath{\OpA}}
\DeclareMathOperator*{\argmin}{arg\,min}
\DeclareMathOperator{\errorfunc}{E}
\DeclareMathOperator{\RegOp}{\mathcal{S}}

\DeclareMathOperator{\loss}{L}

\DeclareMathOperator{\distance}{d}

\DeclareMathOperator{\grad}{\nabla\!}
\DeclareMathOperator{\ProxOp}{\ensuremath{prox}}
\DeclareMathOperator{\esssup}{\ensuremath{ess\,sup}}
\DeclareMathOperator{\FeatureExt}{\ensuremath{\mathcal{R}}}
\newcommand{\stochastic}[1]{\mathsf{#1}}
\newcommand{\stsignal}{\stochastic{\signal}}
\newcommand{\stdata}{\stochastic{\data}}
\DeclareMathOperator{\ProbdistFunctional}{\stochastic{F}}
\DeclareMathOperator{\Expect}{\mathbb{E}}

\newcommand{\ProbabilityMeasure}{\mu}
\newcommand{\signal}{\ensuremath{f}}

\newcommand{\signaltrue}{\signal_{\text{true}}}

\newcommand{\data}{\ensuremath{g}}
\newcommand{\noise}{\delta\data}
\newcommand{\signallearned}[1]{\ForwardOpInvLearned(#1)}
\newcommand{\memory}{s}
\newcommand{\param}{\ensuremath{\Theta}}
\newcommand{\vparam}{\ensuremath{\boldsymbol{\param}}}
\newcommand{\weight}{\ensuremath{w}}

\newcommand{\bias}{\ensuremath{b}}
\newcommand{\primal}{\ensuremath{f}}
\newcommand{\dual}{\ensuremath{h}}
\providecommand{\keywords}[1]{\textbf{\textit{Keywords---}} #1}

\graphicspath{{figures/}}
\usepackage{paperacronyms}
\acuse{MRI,ADMM,GPU,CPU}
\crefname{equation}{}{}
\Crefname{equation}{}{}
\crefname{item}{}{}
\Crefname{item}{}{}
\hypersetup{%
    bookmarks=true,          % show bookmarks bar
    pdfmenubar=true,	        % show Acrobat's menu
    pdffitwindow=false,     % window fit to page when opened
    pdfstartview={FitH},    % fits the width of the page to the window
    colorlinks=true,        % false: boxed links; true: colored links
    linkcolor=red,          % color of internal links
    citecolor=red,          % color of links to bibliography
    filecolor=magenta,      % color of file links
    urlcolor=cyan,          % color of external links
    pdfborder = {0,0,0}
}
\addunit{\pixel}{pixel}
\addunit{\voxel}{voxel}
\addunit{\decibel}{dB}
\addunit{\byte}{B}
\addunit{\hounsfield}{HU}

\title{Solving ill-posed inverse problems using iterative deep neural networks}
\subtitle{General approach illustrated on non-linear tomographic reconstruction \\[1em]
2017-05-22, Version 1.0}
\author{Jonas Adler$^{* \, \dagger}$ \and Ozan \"{O}ktem$^*$}
\institution{%
  {\footnotesize (*) Department of Mathematics, KTH - Royal Institute of Technology} \\ 
  {\footnotesize \phantom{(*)} SE-100 44 Stockholm, Sweden. \texttt{\{jonasadl,ozan\}@kth.se}}  \\ 
  {\footnotesize ($\dagger$) Elekta AB, Box 7593, SE-103 93 Stockholm, Sweden}}

\begin{document}

\maketitle
		
\begin{abstract}
We propose a partially learned approach for the solution of ill posed inverse problems with not necessarily linear forward operators. The method builds on ideas from classical regularization theory and recent advances in deep learning to perform learning while making use of prior information about the inverse problem encoded in the forward operator, noise model and a regularizing functional. The method results in a gradient-like iterative scheme, where the ``gradient'' component is learned using a convolutional network that includes the gradients of the data discrepancy and regularizer as input in each iteration.

We present results of such a partially learned gradient scheme on a non-linear tomographic inversion problem with simulated data from both the Sheep-Logan phantom as well as a head CT. The outcome is compared against \acl{FBP} and \acl{TV} reconstruction and the proposed method provides a \unit{5.4}{\decibel} \acs{PSNR} improvement over the \acl{TV} reconstruction while being significantly faster, giving reconstructions of $\unit{512 \times 512}{\pixel}$ images in about $\unit{0.4}{\second}$ using a single \ac{GPU}.
\end{abstract}

\acresetall
\acuse{MRI,ADMM,GPU,CPU}
	
\section{Introduction}
Inverse problems refer to problems where one seeks to reconstruct parameters characterizing the system under investigation from indirect observations. Such problems arise in several areas of science and engineering. Mathematically, an inverse problem can be formulated as reconstructing (estimating) a signal $\signaltrue \in \RecSpace$ from data $\data \in \DataSpace$ where
\begin{equation}\label{eq:InvProb}
	\data = \ForwardOp(\signaltrue) + \noise.  
\end{equation}  
In the above, $\RecSpace$ and $\DataSpace$ are topological vector spaces, $\ForwardOp \colon \RecSpace \to \DataSpace$ (forward operator) models how a given signal gives rise to data in absence of noise, and $\noise \in \DataSpace$ is a single sample of a $\DataSpace$-valued random variable that represents the noise component of data. 
	
Many inverse problems, such as those arising in imaging, are naturally formulated when both signal and data are functions. In such case, $\RecSpace$ is some Banach/Hilbert space of functions defined on a fixed image domain $\domain \subset \Real^d$ and $\DataSpace$ is likewise a Banach/Hilbert space of functions defined on a fixed data manifold $\datamanifold$. An important remark here relates to the nature of the data manifold. It can be a subset of Euclidean space, but this is not necessarily the case. In fact, x-ray tomographic imaging leads to inverse problems where elements in the data manifold $\datamanifold$ represent lines in $\Real^d$.
	
\subsection{Classical regularization}\label{sec:ClassicalReg}
A common approach in solving an inverse problem of the form in \cref{eq:InvProb} is to minimize the miss-fit against data. For example by minimizing
\begin{equation} 
	\signal \to \LogLikelihood\bigl(\ForwardOp(\signal), \data\bigr)
	\label{eq:datadiscr}
\end{equation}
where $\LogLikelihood \colon \DataSpace \times \DataSpace \to \Real$ is a suitable affine transformation of the data log-likelihood \cite{BeLaZa08}. Then, one may interpret minimizing the above as finding a maximum likelihood solution to \cref{eq:InvProb}.
	
This minimization is a large scale optimization problem that for typical choices of $\ForwardOp$ is ill posed, that is, a solution (if it exists) is unstable with respect to the data $\data$ in the sense that small changes to data results in large changes to a reconstruction. Hence, finding a maximum likelihood solution (there may be several) typically leads to over-fitting against data. 
	
Within classical regularization theory, there are currently three strategies for avoiding over-fitting when solving \cref{eq:InvProb}. One is approximate inverse that is applicable to cases when $\RecSpace$ has a mollifier. The idea is to construct a pseudo-inverse to $\ForwardOp$ using the mollifier \cite{Sc07}. Another approach is iterative regularization, which starts out by considering a fixed point iteration scheme for minimizing \cref{eq:datadiscr}. Over-fitting is avoided by stopping the iterates early, which is a feasible strategy if the iterates are semi-convergent \cite{EnHaNe00,Ha97,BeLaZa08,KaNeSc08}. Finally, we have variational regularization where over-fitting is avoided by introducing a functional $\RegOp \colon \RecSpace \to \Real$ (regularization functional) that encodes a priori information about $\signaltrue$ and penalizes unfeasible solutions \cite{EnHaNe00,ScGrGrHaLe09}. Hence, instead of minimizing only the data discrepancy functional, one now seeks to minimize the regularized objective functional by solving  
\begin{equation}\label{eq:VarReg}
  \min_{\signal \in \RecSpace} \bigl[ \LogLikelihood\bigl(\ForwardOp(\signal), \data\bigr) + \lambda \RegOp(\signal) \bigr]
  \quad\text{for a fixed $\lambda \geq 0$.} 
\end{equation}
In the above, $\lambda$ (regularization parameter) governs the influence of the a priori knowledge encoded by the regularization functional against the need to fit data.
A typical example of variational regularization in imaging is \ac{TV} regularization, which applies to signals that are represented by scalar functions of bounded variation. The corresponding regularization functional is then given as
$\RegOp(\signal) := \Vert \grad \signal \Vert_1$.

\subsection{Machine learning approaches to inverse problems}
\label{sec:machine_learning_inverse}
Machine learning can be seen as algorithms for non-linear function approximation under weak assumptions. Applied to the inverse problem in \cref{eq:InvProb}, it can be phrased as the problem of reconstructing a (non-linear) mapping $\ForwardOpInvLearned \colon \DataSpace \to \RecSpace$ satisfying the following \emph{pseudo-inverse} property:
\[
  \ForwardOpInvLearned(\data) \approx \signaltrue
  \quad\text{whenever data $\data$ is related to $\signaltrue$ as in \cref{eq:InvProb}.}
\]

A key element in machine learning approaches is to parametrize the set of such pseudo-inverse operators by a parameter $\vparam \in \ParamSpace$. The "learning" part refers to choosing an "optimal" parameter given some \emph{training data}, where the concept of optimality is typically quantified through a loss functional that measures the quality of a learned pseudo-inverse $\ForwardOpInvLearned$. 

The manner in which the loss functional is specified depends on the type of training data, and here we separate between supervised and unsupervised learning. These two approaches are \emph{fundamentally} different and this article focuses on the supervised learning case since it is the problem with the most structure and we expect learning to give larger improvements over traditional methods.

\subsubsection{Supervised learning}
In supervised learning, training data are independent identically distributed realizations of a $(\DataSpace\times \RecSpace)$--valued random variable $(\stdata,\stsignal)$ with a known probability density $\ProbabilityMeasure$.
Estimating $\vparam \in \ParamSpace$ from training data can be formulated as minimizing a loss functional $\vparam \mapsto \loss(\vparam)$ that frequently has the following structure: 
\begin{equation}\label{eq:LossAbstract}
  \loss(\vparam) := \ProbdistFunctional_{\ProbabilityMeasure} \Bigl[ \distance\bigl(\signallearned{\stdata}, \stsignal \bigr) \Bigr].
\end{equation}
In the above, $\ForwardOpInvLearned \colon \DataSpace \to \RecSpace$ is the pseudo-inverse that is given by $\vparam \in \ParamSpace$, 
$\distance \colon \RecSpace \times \RecSpace \to \Real$ is a \emph{"distance" function} quantifying the quality of a specific reconstruction, and $\ProbdistFunctional_{\ProbabilityMeasure}$ maps real-valued random variables on $\DataSpace\times \RecSpace$ to real numbers. 

A common choice is to use the expected loss w.r.t. the squared distance:
\begin{equation}\label{eq:LossStandard}
  \loss(\vparam) :=
    \Expect_{\ProbabilityMeasure} 
    \Bigl[
	\bigl\Vert \signallearned{\stdata} - \stsignal \bigr\Vert_\RecSpace^2
    \Bigr].
\end{equation}
One may also consider other loss functionals and the method we suggest can easily be adapted to these. As an example, a very conservative reconstruction method would use a loss functional given by the supremum of the $\infty$-norm:
\[
  \loss(\vparam) :=
    \esssup_{\ProbabilityMeasure} 
    \Bigl[
	\bigl\Vert \signallearned{\stdata} - \stsignal \bigr\Vert_\infty
    \Bigr].
\]

\subsubsection{Unsupervised learning}
In unsupervised learning there is no access to input-output pairs as in the supervised learning setting. Instead, here training data are given as elements in $\DataSpace$. The corresponding mathematical setting is to consider these as independent identically distributed realizations of a $\DataSpace$--valued random variable $\stdata$ with a known probability density $\ProbabilityMeasure$. A natural choice for a loss function is now to quantify how well the learned reconstruction matches the regularized data discrepancy, i.e.,
\[
  \loss(\vparam) := \Expect_{\ProbabilityMeasure} \Bigl[ \LogLikelihood\Bigl(\ForwardOp\bigl(\signallearned{\stdata}\bigr), \stdata\Bigr) + \RegOp\bigl(\signallearned{\stdata} \bigr) \Bigr].
\]
The above can be interpreted as learning an optimizer for the variational regularization in \cref{eq:VarReg}.
	
\subsection{Survey of the field}
From the data science perspective, the classical regularization approaches outlined in \cref{sec:ClassicalReg} are all examples of knowledge-driven modelling whereas machine learning is usually categorized as data-driven modelling. These two modelling paradigms can be combined in different way in order to solve \cref{eq:InvProb}.
	
\paragraph{Fully learned reconstruction}
Approaching the inverse problem directly with machine learning amounts to learn $\ForwardOpInvLearned \colon \DataSpace \to \RecSpace$ from data such that it approximates an inverse of $\ForwardOp$ in \cref{eq:InvProb}. An example of such an approach for solving small scale tomographic reconstruction problems is given in \cite{PaGiKaLoMa04, ArMaTsSt12}.
	
An obvious disadvantage with such fully learned approaches is that the result is likely to depend on the data manifolds, so training data needs to be rich enough to account for all various data manifolds that one is likely to encounter. Furthermore, training data also needs to be rich enough to allow the learning scheme to learn the structure in $\ForwardOp$, which is given by the physics laws governing the formation of data from a signal. Finally, in many applications the adequate digitalizations of the signal and data often requires very high dimensional arrays.
	
The above considerations imply that the parameter space $\ParamSpace$ used for parametrizing possible inverse operators has to be very high dimensional in a fully learned approach. Therefore, the idea of learning $\ForwardOpInvLearned$ from data without using any knowledge of the physics quickly or the data manifold quickly becomes in-feasible due to lack of training data. 
	
\paragraph{Sequential data and knowledge driven reconstruction}
The idea here is to separate the learned components from a part that encodes some knowledge about the structure of $\ForwardOp$ and the data manifold. Formalizing this, we assume 
\begin{equation}\label{eq:SepInvLearn}
  \ForwardOpInvLearned = \OpB_{\vparam} \circ \OpA \circ \OpC_{\vparam}
\end{equation}
where $\OpA  \colon \DataSpace \to \RecSpace$ is a known component that encodes knowledge about the structure of the forward operator $\ForwardOp$ whereas the operators $\OpB_{\vparam} \colon \RecSpace \to \RecSpace$ and $\OpC_{\vparam} \colon \DataSpace \to \DataSpace$ are the learned components.
	
An important special case is $\OpC_{\vparam} = \IdentityOp$, which entirely separates the computation of $\OpA$ from the learning of $\OpB_{\vparam}$, i.e., the original inverse problem in \cref{eq:InvProb} is recast as learning $\OpB_{\vparam} \colon \RecSpace \to \RecSpace$. This significantly simplifies the implementation since the often demanding task of computing $\data \mapsto \OpA(\data)$ can be separated from the learning software. It also ensures the data manifold is not explicitly part of the learning. One may furthermore exploit additional structure, like locality and/or invariance of the operator $\OpA \circ \ForwardOp$.  
A key step for such a sequential data and knowledge driven reconstruction scheme is to have candidates for $\OpA$ and one natural option is to let it be some pseudo-inverse. As an example, in tomographic applications it can be given by the back-projection or the filtered back-projection operator. Next, when learning $\OpB_{\vparam}$ from data, one may use approaches that build on the corpus of knowledge that exists for denoising signals in $\RecSpace$. An example demonstrating this approach for tomographic reconstruction is given in \cite{JiMcFrUn16,HaWuPoMa17}.

On the other hand, for ill posed inverse problems some information is irreversibly lost when making the assumption in \cref{eq:SepInvLearn} since the learned operators $\OpC_{\vparam}$ and $\OpB_{\vparam}$ cannot recover information that is lost by using a pseudo-inverse $\OpA$. 
To alleviate this problem, for linear forward operators we can consider choosing $\OpA$ as the adjoint of the forward operator, i.e., $\OpA := \ForwardOp^* \colon \DataSpace \to \RecSpace$. This can be seen as learning to solve the normal equations for \cref{eq:InvProb} since 
\[
  \signaltrue \approx \OpB_{\vparam} \bigl(\ForwardOp^* \bigl(\data \bigr)\bigr) 
  \implies
  \signaltrue \approx \OpB_{\vparam} \bigl(\ForwardOp^* \bigl(\ForwardOp(\signaltrue)\bigr)\bigr)
  \implies
  \OpB_{\vparam} \approx (\ForwardOp^* \circ \ForwardOp)^{-1}.
\]
Nonetheless, solving the normal equations for ill-posed problems is often more ill-posed than the original inverse problem, so such a learning procedure would need to include some kind of regularization.

We conclude by pointing to examples where the operator $\OpC_{\vparam}  \colon \DataSpace \to \DataSpace$ is learned. One such case is \cite{WuGhChMa16} where tomographic reconstruction is performed by learning $\OpC_{\vparam}$. Another similar, but more advanced, example is \cite{PeBa13} where the operator $\OpA$ is given by several \ac{FBP} operators and the learned operator $\OpC_{\vparam}$ is given by the filter coefficients.
	
\paragraph{Learning an iterative reconstruction scheme}
The previous method for combining knowledge and data driven reconstruction is ultimately limited by what knowledge $\OpA \colon \DataSpace \to \RecSpace$ manages to capture about the inverse of $\ForwardOp \colon \RecSpace \to \DataSpace$. 

To address this limitation one may formulate a parametrized optimization problem for solving \cref{eq:InvProb}, typically of the form \cref{eq:VarReg}, and then assign suitable values to these parameters by learning. This procedure can be re-formulated as a bi-level optimization scheme whose mathematical properties (like existence) can be analysed in a functional analytic setting \cite{ReSc13,CaCaReScVa15,CaReSc16,ReScVa16,ReScVa17}. Obviously any implementation of an iterative scheme for solving the aforementioned parametrized optimization problem will terminate after a finite number if iterates. Hence, the outcome will not only depend on the formulation of the parametrized optimization problem, but also on the solution scheme one chooses to use. Hence, the above bi-level optimization scheme by itself does not uniquely determine a reconstruction operator for \cref{eq:InvProb}. 

A natural extension of the bi-level optimization scheme above is to include the solution scheme in the learning. One such approach was given in \cite{AnDeCoHoPf16}. Here, an iterative stochastic gradient method is learned from data consisting of optimization problems, each associated with a deep learning problem. The output is thus a trained stochastic gradient method that can be used to train other deep neural networks. This overall "learning to learn" approach can in principle be extended to other use cases by merely changing the underlying data type. An example of such an extension to solving inverse problems is \cite{YaSuLiXu16}, which learns an \ac{ADMM}-like scheme for \ac{MRI} reconstruction. Another is \cite{ChLiPoViSa17}, which learns a "proximal" in an \ac{ADMM}-like scheme for various image restoration problems. Finally, \cite{PuWe17} considers solving finite dimensional linear inverse problems typically arising in image restoration. The idea is to learn over a broader class of schemes instead of restrict attention to a specific type of scheme, like \ac{ADMM} above.
	
\subsection{Contribution and overview of paper}
This paper generalizes the ideas in \cite{PuWe17} in many directions. First, we consider solving (possibly) non-linear inverse problems in a functional analytic setting. Next, we consider the issue of proper initialization and include further prior knowledge by allowing a regularizer. Finally, we also provide a generic and scalable open source implementation\footnote{
\href{https://github.com/adler-j/learned\_gradient\_tomography}{https://github.com/adler-j/learned\_gradient\_tomography}} of our method based on \ac{ODL} \cite{AdKoOk17} that can be applied to a wide range of realistic inverse problems. We also provide the trained parameter $\vparam$ used for generating the results shown in the article. To show that the approach can handle (non-linear) forward operators in large scale inverse problems, we consider tomographic reconstruction with a non-linear ray transform inversion.
	
\Cref{sec:GeneralAlg} derives a partially learned gradient decent scheme for solving \cref{eq:InvProb} in the functional analytic setting. This section also introduces the deep convolutional network that is used later for tomographic reconstruction.
\Cref{sec:Impl} describes the implementation of the partially learned gradient decent scheme in \cref{sec:GeneralAlg} and software components used for computations. \Cref{sec:Eval} tests the performance of the  partially learned gradient decent scheme on tomographic inverse problems. The paper concludes with a discussion in \cref{sec:Disc} and a summary of future work and conclusions is given in \cref{sec:Conc}.

\section{Solving inverse problems by learned gradient descent}\label{sec:GeneralAlg}
We begin by proving a heuristic motivation that comes from comparing two natural considerations involving gradient descent schemes associated with solving \cref{eq:InvProb}. This results in an initial scheme for partially learned gradient descent given in \cref{alg:main_alg1}, which is then extended by adding persistent memory and resulting in the scheme in \cref{alg:main_alg}. Next is a description of how to integrate deep learning with \cref{alg:main_alg}, resulting in the final scheme given in \cref{alg:learnedfull}.

\subsection{Motivation}
The starting point in learning an iterative scheme that combines elements from (deep) machine learning and classical regularization theory is to consider the \emph{error functional} $\errorfunc \colon \RecSpace \to \Real$ defined as 
\[  \errorfunc(\signal) := \distance(\signal,\signaltrue)  \]
where $\distance \colon \RecSpace \times \RecSpace \to \Real$ is the distance functional that appears in the definition of the loss functional \cref{eq:LossAbstract} used for training. It measures how good well $\signal$ approximates $\signaltrue$, so one natural error functional corresponding to \cref{eq:LossStandard} is
\[ 
  \errorfunc(\signal) = \Vert \signal - \signaltrue \Vert_{\RecSpace}^2.
\] 

Ideally, solving \eqref{eq:InvProb} would be based on minimizing the error functional, which for obvious reasons is not possible. Hence, we need to use a substitute. In variational regularization theory, such a substitute is given by the regularized objective functional in \cref{eq:VarReg}. Much of regularization theory aims at choosing the objective functional in \ref{eq:VarReg} so that the regularized solution approximates the true signal to be recovered:
\begin{equation}\label{eq:ApproxEqVar}
  \argmin_{\signal\in \RecSpace} \errorfunc(\signal) \approx 
  \argmin_{\signal\in \RecSpace} \Bigl[  \LogLikelihood\bigl(\ForwardOp(\signal), \data\bigr) + \lambda \RegOp(\signal) \Bigr].
\end{equation}
Assume next that the objective functional in the right hand side of \cref{eq:ApproxEqVar} is Fr\'echet differentiable and (strictly) convex. Then, a simple gradient descent scheme could be used to find a minimum:
\begin{equation}\label{eg:GDRHS}
  \signal_{i} := \signal_{i-1} - \sigma \Bigl( \grad \big[ \LogLikelihood\bigl(\ForwardOp(\cdot), \data\bigr) \big](\signal_{i-1}) + \lambda [\grad\RegOp](\signal_{i-1}) \Bigr)
\end{equation}
where, assuming a differentiable likelihood and forward operator, we note that
\[
	\grad \big[ \LogLikelihood\bigl(\ForwardOp(\cdot), \data\bigr) \bigr](\signal)
	=
	[\partial \ForwardOp](\signal)^*
	\Bigl(
		[\grad \LogLikelihood(\cdot, \data)]\bigl(\ForwardOp(\signal) \bigr)
	\Bigr)
	\quad
	\text{for any $\signal \in \RecSpace$.}
\]
Likewise, considering the left hand side in the same way, a differentiable convex error functional would allow one to use a corresponding scheme for finding a minimum:
\begin{equation}\label{eg:GDLHS}
  \signal_{j+1} := \signal_j - \sigma [\grad\errorfunc](\signal_j).
\end{equation}
Since the gradient mapping $\grad\errorfunc \colon \RecSpace \to \RecSpace$ in \cref{eg:GDLHS} requires knowledge about the true signal, it is natural to try to learn it from training data while utilizing knowledge about the gradient mappings 
$\grad \big[ \LogLikelihood\bigl(\ForwardOp(\cdot), \data\bigr) \bigr], \grad\RegOp \colon \RecSpace \to \RecSpace$. For this purpose, we introduce the \emph{(learned) updating operator}    
$\Lambda_{\vparam} \colon \RecSpace \times \RecSpace \times \RecSpace \to \RecSpace$
that, given an appropriately selected (learned) parameter $\vparam\in \ParamSpace$, should satisfy
\[  
  \Lambda_{\vparam}\Bigl(
    \signal, \grad \big[ \LogLikelihood\bigl(\ForwardOp(\cdot), \data\bigr) \bigr](\signal), \lambda \grad \RegOp(\signal)
  \Bigr)
  \approx
  \grad \errorfunc(\signal).
\]
These considerations suggests a partially learned gradient descent scheme specified as  in \cref{alg:main_alg1}.
\begin{algorithm}
\caption{Initial partially learned gradient descent}\label{alg:main_alg1}
\begin{algorithmic}[1]
  \State Select an initial guess $\signal_0$	
  \For{$i = 1, \dots, I$}
    \State $\Delta \signal_i \gets -\sigma
      \Lambda_{\vparam}\Bigl(
         \signal_{i-1}, \grad \big[ \LogLikelihood\bigl(\ForwardOp(\cdot), \data\bigr) \bigr](\signal_{i-i}), \lambda \grad \RegOp(\signal_{i-1})
      \Bigr)$
    \State $\signal_i \gets \signal_{i-1} + \Delta \signal_i$
  \EndFor
  \State $\ForwardOpInvLearned(g) \gets \signal_I$
\end{algorithmic}
\end{algorithm}

\Cref{alg:main_alg1} suffers from several unnecessary shortcomings that are easily addressed. 
The regularization parameter $\lambda$ and the step length $\sigma$ have to be explicitly chosen, a task that is known to be troublesome in practical applications (\cref{sec:UseCases}). One may instead make these part of $\vparam$ and thereby learn them from training data. 
Next, the convergence rate of gradient descent schemes can be accelerated by using information from previous iterates (memory) as in quasi-Newton schemes\cite{LBFGS}. For this purpose we introduce persistent memory $\memory \in \RecSpace^M$ that allows \cref{alg:main_alg1} to use information from earlier iterates. The learned updating operator now becomes a mapping 
\begin{equation}\label{eq:NewLearnedUpdate} 
\Lambda_{\vparam} \colon \RecSpace^M \times \RecSpace \times \RecSpace \times \RecSpace  \to \RecSpace^M \times \RecSpace.  
\end{equation} 
Finally, one often also has the possibility to select the initial iterate $\signal_0$ using some suitable pseudo-inverse $\ForwardOpPseudoInv \colon \DataSpace \to \RecSpace$. 
Considering these modifications results in the partially learned gradient descent scheme listed in \cref{alg:main_alg}.
\begin{algorithm}
\caption{Partially learned gradient descent}\label{alg:main_alg}
\begin{algorithmic}[1]
  \State $\signal_0 \gets \ForwardOpPseudoInv(\data)$.	
  \State Initialize ``memory'' $\memory_0 \in \RecSpace^M$.
  \For{$i = 1, \dots, I$}
    \State $(\memory_i, \Delta \signal_i) \gets 
    	\Lambda_{\vparam}\Bigl(
  	  \memory_{i-1}, \signal_{i-1}, \grad \big[ \LogLikelihood\bigl(\ForwardOp(\cdot), \data\bigr) \big](\signal_{i-1}), \grad \RegOp(\signal_{i-1})
	\Bigr)	$
    \State $\signal_i \gets \signal_{i-1} + \Delta \signal_i$
  \EndFor
  \State $\ForwardOpInvLearned(g) \gets \signal_I$
\end{algorithmic}
\end{algorithm}

\subsection{Parametrizing the learned updating operators}
\label{sec:hyperparam}
The goal here is to specify the class of learned updating operators that are parametrized by $\vparam \in \ParamSpace$.
Following the paradigm in (deep) neural networks, we start by defining a family of affine operators 
\begin{equation}\label{eq:AffineOp}
  \AffineOp_{\weight_n,\bias_n} \colon \RecSpace^{c_{n-1}} \to \RecSpace^{c_{n}}
  \quad\text{for $n=0, \ldots, N$,}
\end{equation}
parametrized by linear mappings $w_n \colon \RecSpace^{c_{n-1}} \to \RecSpace^{c_{n}}$ (weights) and $b_n \in \RecSpace^{c_n}$ (biases). Here, $N$ is usually referred to as  the depth of the neural network that will eventually define the learned updating operator and $c_n$ is the number of channels in the $n$:th layer. 
Next, we introduce a family of non-linear operators  
\begin{equation}\label{eq:NonLinOp}
\NonLinOp_{n} \colon \RecSpace^{c_{n}} \to \RecSpace^{c_{n}}
\end{equation}  
that are given by point-wise application of a fixed non-linear scalar function, henceforth called the response function. 

By chaining compositions, we now define a parametrized family of learned updating operators as
\[
\Lambda_{\vparam} := (\NonLinOp_{N} \circ \AffineOp_{\weight_N,\bias_N}) \circ \cdots \circ (\NonLinOp_{1} \circ \AffineOp_{\weight_1,\bias_1})
\]
with $\vparam := \bigl((w_N, b_N), \dots, (w_1, b_1) \bigr)$. In order to match the  domain and range of the operator in \cref{eq:NewLearnedUpdate}, we need to assume that $c_0 = M + 3$ and $c_N = M + 1$. 
	
Such parametrized operators are used in machine learning applications for two primary reasons: computability and descriptive power. In order to learn the parameters $\vparam$ from training data, a (stochastic) gradient descent method is typically applied in which case the derivative $\partial \Lambda_{\vparam} / \partial \vparam$ and its adjoint $[\partial \Lambda_{\vparam} / \partial \vparam]^*$ needs to be repeatedly computed and here one may use the chain rule. This becomes particularity easy to perform in a  computationally feasible manner for learned updating operators of this form. 
Furthermore, introducing the non-linear component in \cref{eq:NonLinOp} allows the learned operator to approximate a large set of non-linear operators \cite{Gu16}. 

\paragraph{Choice of affine and non-linear operator families}
Our next step is to further narrow down the generative models for the operator families $\AffineOp_{\weight_i, \bias_i}$ and $\NonLinOp_i$.  
We start by writing the affine operator in \cref{eq:AffineOp} as
\[ \AffineOp_{\weight_n,\bias_n} = 
     ( \AffineOp^{1}_{\weight_n,\bias_n}, 
       \ldots, 
       \AffineOp^{c_n}_{\weight_n,\bias_n}
     )
\]
where the components  
\[ 
\AffineOp^{l}_{\weight_n,\bias_n} \colon \RecSpace^{c_{n-1}} \to \RecSpace
\quad\text{for $l=1,\ldots, c_n$}
\]
represent the affine transformation for the $l$:th channel in the $n$:th layer. 

Furthermore, for many inverse problems it is sufficient to assume that $\AffineOp^{l}_{\weight_n,\bias_n}$ is translation invariant, which implies that 
\[ 
 \AffineOp^{l}_{\weight_n,\bias_n}(\signal_1,\ldots,\signal_{c_{n-1}}) 
 = \bias_n^l + \sum_{j=1}^{c_{n-1}} \weight_{n}^{j,l} \ast \signal_j 
\]
where $\bias_n^l \in \Real$ represents the bias and $\weight_n$ is given as a ``matrix'' of convolution kernels $\weight_{n}^{j,l} \in \RecSpace$. Hence, our parameter space becomes 
\[  
\ParamSpace = (\RecSpace^{c_N \times c_{N-1}} \times \Real^{c_N}) \times \ldots \times (\RecSpace^{c_1 \times c_{0}} \times \Real^{c_1}).
\]

Finally, the non-linear response functions $\NonLinOp_i$ in \cref{eq:NonLinOp} can be chosen in different ways and we will be using the \ac{ReLU} \cite{relu}
\[
\text{relu}(x) = \begin{cases}
x & \text{if } x > 0 \\
0 & \text{else.}
\end{cases}
\]

\subsection{The partially learned gradient descent algorithm}
A number of hyper-parameters needs to be chosen prior to learning. These are the number of layers $N \in \Natural$, the number of channels $c_1, \ldots, c_{N-1} \in \Natural$ in each layer, the number of iterations $I$, and the size of the memory $M$.  

In the examples shown in in \cref{sec:Eval}, we let the weights $\weight_i$ be represented by \unit{3 \times 3}{\pixel} convolutions and we used $N=3$ layers. The number of convolutions in each layer was selected as $m_1 = 32$ and $m_2 = 32$. We selected the number of iterations to be $I=10$ and the amount of memory to be $M=5$. Such low numbers were selected in order to reduce the space of allowed parameters which in turn should help reduce over-fitting. All parameters were selected by simple trial and error and it is likely that a more sophisticated set-up would give better results or be better suited for a particular application. For further details, see the supplemental source code.

Once the hyper-parameters are chosen, one may now fully specify the partially learned gradient descent, which is done in \cref{alg:learnedfull}. The scheme learns the updating operator by learning the scalars $\bias_n^l \in \Real$ and the convolution kernels (functions) $\weight_{n}^{j,l}$ from training data. The resulting learned updating operator can then be used to solve the inverse problem in \cref{eq:InvProb}.
\begin{algorithm}
	\caption{Partially learned gradient descent}\label{alg:learnedfull}
	\begin{algorithmic}[1]
		\State $\memory_0 \gets 0$
		\State $\signal_0 \gets \ForwardOpPseudoInv(\data)$
		\For{$i = 1, \dots, I$}
		\State $u_i^1 \gets \Bigl(\signal_{i-1}, \memory_{i-1},  \grad \big[ \LogLikelihood\bigl(\ForwardOp(\cdot), \data\bigr) \big](\signal_{i-1}), \grad \RegOp(\signal_{i-1}) \Bigr)$
		\State $u_i^2 \gets \text{relu}\bigl( \AffineOp_{\weight_1,\bias_1}(u_i^1) \bigr)$ 		
		\State $u_i^3 \gets \text{relu}\bigl( \AffineOp_{\weight_2,\bias_2}(u_i^2) \bigr)$
		\State $( u_i^4, \Delta \signal_i ) \gets \AffineOp_{\weight_3,\bias_3}(u_i^3)$
		\State $\memory_i \gets \text{relu}(u_i^4)$
		\State $\signal_i \gets \signal_{i-1} + \Delta \signal_i$
		\EndFor
		\State $\ForwardOpInvLearned(g) \gets \signal_I$
	\end{algorithmic}
\end{algorithm}

\section{Implementation and evaluation}\label{sec:Eval}
The algorithm was tested on the two-dimensional computed tomography problem. The signal is in this case real valued functions defined on a domain in $\Real^2$ representing images and $\RecSpace$ is a suitable vector space of such functions.
The corresponding forward operator is expressible in terms of the ray transform $\RadonTransform : \RecSpace \to \DataSpace$, which integrates the signal over a set of lines $\datamanifold$ given by the acquisition geometry. Hence, elements in $\DataSpace$ are functions on lines
\[
	\RadonTransform(\signal)(\ell) = \int_{\ell} f(x) dx \quad \text{ for $\ell \in \datamanifold$.}
\]

\begin{figure}
	\centering	
	\begin{subfigure}[t]{.32\linewidth}
		\includegraphics[width=\linewidth, trim={21mm 16mm 27mm 6mm}, clip]{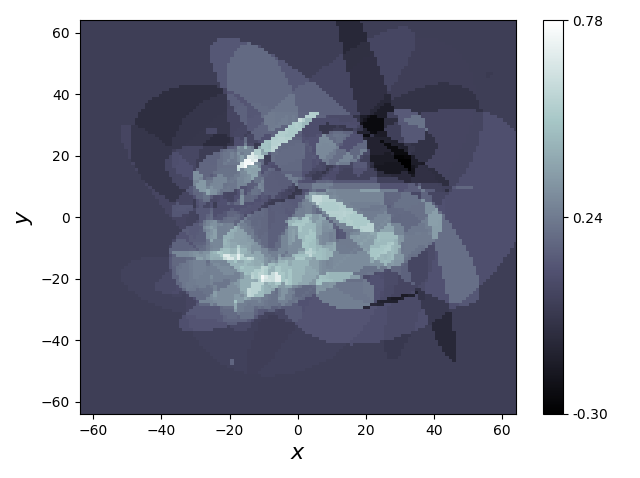}
		\label{fig:random_ellipse_phantom}
	\end{subfigure}
	\begin{subfigure}[t]{.32\linewidth}
		\includegraphics[width=\linewidth, trim={21mm 16mm 27mm 6mm}, clip]{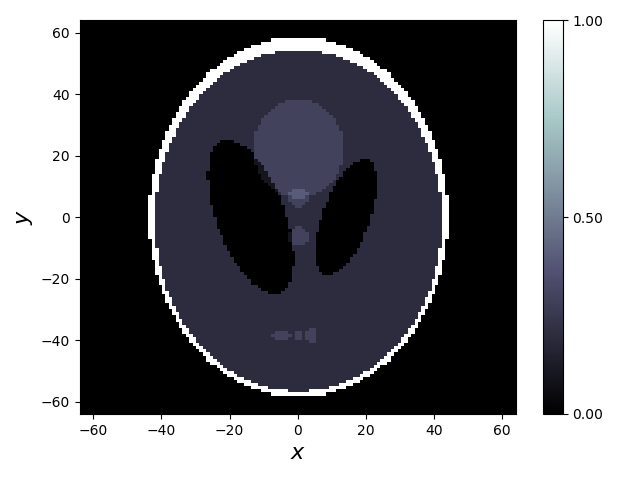}
		\label{fig:shepp_logan_phantom}
	\end{subfigure}
	\begin{subfigure}[t]{.32\linewidth}
		\includegraphics[width=\linewidth, trim={22mm 16mm 27mm 6mm}, clip]{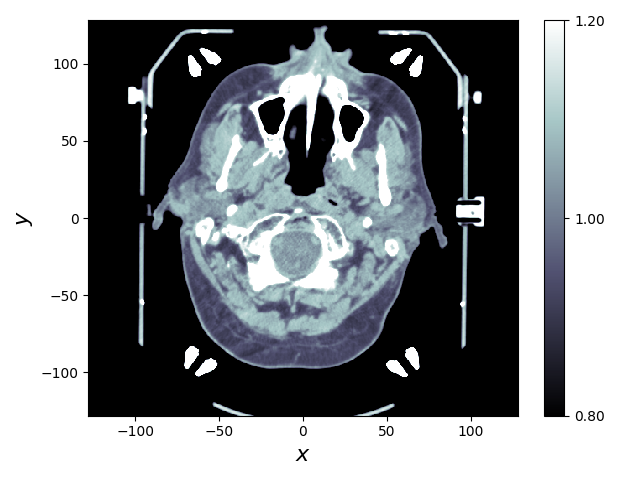}
		\label{fig:head_phantom}
	\end{subfigure}
        \\[1em]
	\begin{subfigure}[t]{.32\linewidth}
		\includegraphics[width=\linewidth, trim={21mm 16mm 27mm 6mm}, clip]{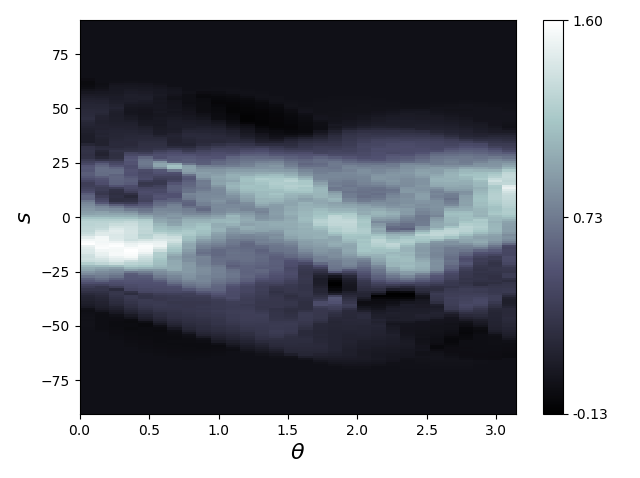}
		\label{fig:random_ellipse_data}
	\end{subfigure}
	\begin{subfigure}[t]{.32\linewidth}
		\includegraphics[width=\linewidth, trim={21mm 16mm 27mm 6mm}, clip]{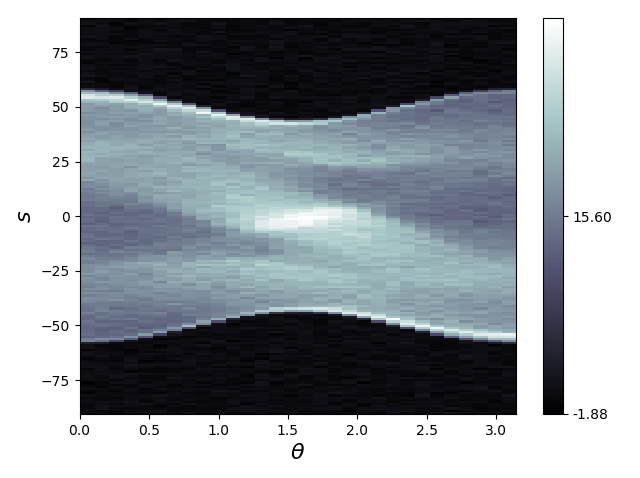}
		\label{fig:shepp_logan_data}
	\end{subfigure}
	\begin{subfigure}[t]{.32\linewidth}
		\includegraphics[width=\linewidth, trim={22mm 16mm 27mm 6mm}, clip]{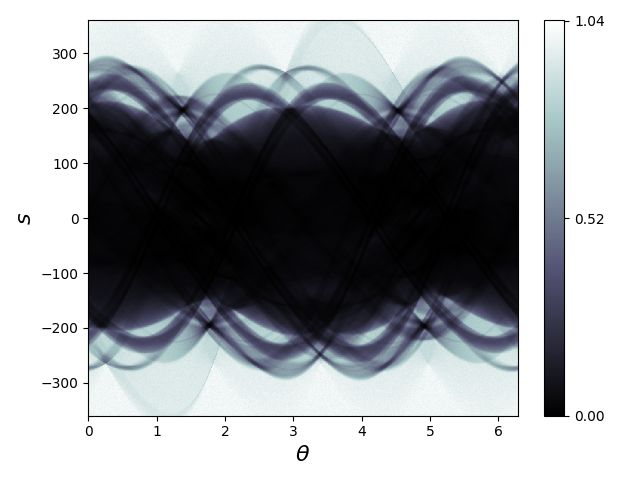}
		\label{fig:head_data}
	\end{subfigure}
        \\[1em]	
	\begin{subfigure}[t]{.32\linewidth}
		\includegraphics[width=\linewidth, trim={21mm 16mm 27mm 6mm}, clip]{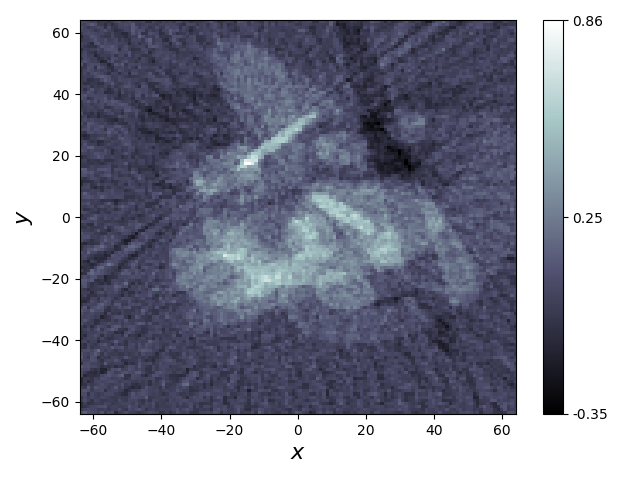}
		\label{fig:random_ellipse_fbp}
	\end{subfigure}
	\begin{subfigure}[t]{.32\linewidth}
		\includegraphics[width=\linewidth, trim={21mm 16mm 27mm 6mm}, clip]{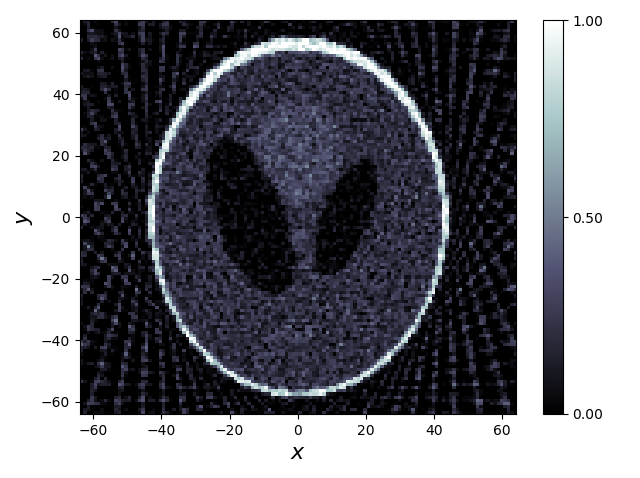}
		\label{fig:shepp_logan_fbp}
	\end{subfigure}
	\begin{subfigure}[t]{.32\linewidth}
		\includegraphics[width=\linewidth, trim={22mm 16mm 27mm 6mm}, clip]{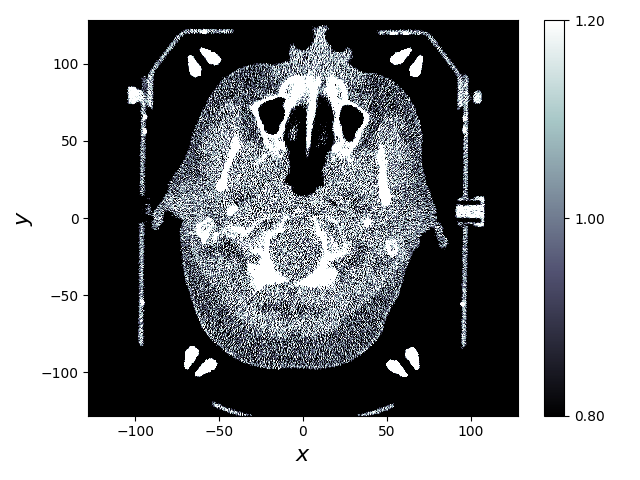}
		\label{fig:head_initial}
	\end{subfigure}
	\caption{Examples from training data. Top row shows the phantoms, middle row the simulated tomographic data, 
	  and the bottom row is the initial guess obtained using \ac{FBP}. The left column is the random ellipses, the middle 
	  column is the Shepp-Logan phantom, and the right column is slice of head phantom shown with 
	  window [-200, 200] \hounsfield.}
	\label{fig:example_data}
\end{figure}

As training data we consider \ac{CT} simulations from two particular types of phantoms with different forward operators and noise models:
\begin{description}
\item[Ellipses:]
Training data is randomly generated ellipses on a \unit{128 \times 128}{\pixel} domain. The projection geometry was selected as a sparse 30 view parallel beam geometry with 5\% additive Gaussian noise added to the projections. In this case, the log-likelihood was selected as the squared $L^2$ norm $\LogLikelihood\bigl(\cdot, \data\bigr) := \frac{1}{2} \bigl\Vert \cdot - \data \bigr\Vert_\DataSpace^2$ which implies 
\[
     \grad \big[ \LogLikelihood\bigl(\RadonTransform(\cdot), \data\bigr) \big](\signal) = \RadonTransform^*\bigl(\RadonTransform(\signal) - \data\bigr)
\]
The phantoms were generated "on the fly", giving an effectively infinite dataset.
\item[Heads:]
The training is simulated projections of \unit{512 \times 512}{\pixel}, \unit{256 \times 256}{\milli\meter} slices of \ac{CT} scans of human heads as provided by Elekta (Elekta AB, Stockholm, Sweden). The acquisition geometry defining the data manifold was selected as a fan beam geometry with source-axis distance of \unit{500}{\milli\meter}, source-detector distance \unit{1000}{\milli\meter}, \unit{1000}{\pixel}, and 1000 angles.

Here, in order to get a accurate noise model we used a non-linear forward operator given by
\[
  \ForwardOp(\signal)(\ell) = \lambda \exp\bigl(-\mu \RadonTransform(\signal)(\ell)\bigr)
\]
where $\lambda\in \Real^+$ is the mean number of photons per \pixel, taken to be $10\,000$, and $\mu \in \Real^+$ is the linear attenuation coefficient which was taken to be that of water ($\approx \unit{0.2}{\centi\reciprocal\meter}$). Poisson noise was added to the projections, and given $10\,000$ photons per \pixel, which corresponds to a low dose scan. For this type of noise, the log-likelihood is given by the Kullback-Leibler divergence and the data discrepancy becomes
\[
  \LogLikelihood\bigl(\ForwardOp(\signal), \data\bigr) := 
		\int_{\datamanifold} \Biggl( \ForwardOp(\signal)(\ell) + \data(\ell) \log\biggl( \frac{\data(\ell)}{\ForwardOp(\signal)(\ell)}\biggr) \Biggr) d\ell
\]
which implies that 
\[
  \grad \big[ \LogLikelihood\bigl(\ForwardOp(\cdot), \data\bigr) \big](\signal) = 
		[\partial\ForwardOp(\signal)]^*\biggl(1.0 - \frac{\data}{\ForwardOp(\signal)} \biggr).
\]
In the above, the adjoint of the derivative of the forward operator applied in a perturbation $\delta \data \in \DataSpace$ is given by
\[
  [\partial\ForwardOp(\signal)]^*(\delta \data) = \RadonTransform^* 
	 \biggl( 
		 - \mu \underbrace{\lambda \exp\bigl(-\mu\RadonTransform(\signal)(\cdot)\bigr)}_{\ForwardOp(\signal)(\cdot)}
		 \delta \data(\cdot)
	 \biggr),
\]
which after some simplifications gives the following expression for the gradient:
\[
  \grad \big[ \LogLikelihood\bigl(\ForwardOp(\cdot), \data\bigr) \big](\signal)
  =
  -\mu \RadonTransform^* 
	\bigl(
		\ForwardOp(\signal) - \data
	\bigr)
	\quad\text{for $\signal \in \RecSpace$.}
\]
The training used 500 \ac{CT} scans with a total of 41\,000 slices.
\end{description}
For both cases the regularizer was selected as the Dirichlet energy, e.g., 
\[
\RegOp(\signal) := \frac{1}{2} \Vert \grad \signal \Vert_2^2 \implies 
\grad \RegOp (\signal) = \grad^{\,*} (\grad \signal)
\]
which is intended to assist the solver in finding edges since this is a typical feature of interest. See \cref{fig:example_data} for examples of the data used for training and validation.

\subsection{Implementation}\label{sec:Impl}
The methods described above were implemented in Python using \ac{ODL} \cite{AdKoOk17} and Tensorflow \cite{AbAgBaBrCh15}. All operator-related components, such as the forward operator $\ForwardOp$, were implemented in \ac{ODL}, and these were then converted into Tensorflow layers using the \texttt{as\_tensorflow\_layer} functionality of \ac{ODL}. The neural network layers and training were implemented using Tensorflow.

The implementation utilizes abstract \ac{ODL} structures for representing functional analytic notions and is therefore generic, yet easily adaptable to other inverse problems.
We used the \ac{ODL} operator \texttt{RayTransform} in order to evaluate the ray transform and its adjoint using the \ac{GPU} accelerated \texttt{'astra\_gpu'} backend \cite{AaPaCaJaBl16}. The pseudo-inverse $\ForwardOpPseudoInv$ was given by the filtered back-projection algorithm implemented in \ac{ODL} as \texttt{fbp\_op} with no additional smoothing filter.

We emphasize that the functional analytic formulation of \cref{alg:learnedfull} is critical to handle problems of this scale. As an example, storing the ray transform used for the heads dataset as a sparse matrix of floating point numbers would require about $\unit{10}{\giga\byte}$ of \ac{GPU} memory.

\paragraph{Training}

We trained the parameters $\vparam$ using the \texttt{RMSPropOptimizer} optimizer in Tensorflow. We initially used $10^5$ batches on the ellipses problem, where each batch contained 20 tomography problems with a learning rate starting at $10^{-3}$ and decayed according to the inverse of the iteration number down to about $10^{-5}$. This training took four days on a workstation with a single Nvidia GTX Titan \ac{GPU}.
These parameters were then used as an initial guess for the heads problem, once again trained according to the same scheme but with the learning rate starting at $10^{-5}$ and decreased to about $10^{-7}$ and with each batch containing only one tomographic problems due to memory limitations of the current implementation. This training took four days on the aforementioned hardware.

\paragraph{Comparison} 

We compare the performance of the partially learned iterative algorithm against the \ac{FBP} algorithm and \ac{TV} regularization. The \ac{FBP} reconstructions were performed using a Hann filter with bandwidth selected to maximize \ac{PSNR}. The data discrepancy in the \ac{TV} regularization matched the one used in the partially learned algorithm and the regularization parameter was selected to maximize \ac{PSNR}.

We solve the \ac{TV} regularized problem without smoothing using the the generic \ac{ODL} implementation of the non-linear Chambolle-Pock primal dual optimization method \cite{Va14}. This is needed since the forward operator in the heads dataset is non-linear. We used 1000 iterations at which point the objective function was stationary. For the ellipses dataset, the evaluation was performed on the modified Shepp-Logan phantom, while on the heads dataset a slice through the nasal region was used.

\subsection{Results}
\begin{figure}[h]
	\centering	
	\begin{subfigure}[t]{.49\linewidth}
		\includegraphics[width=\linewidth, trim={21mm 16mm 32mm 6mm}, clip]{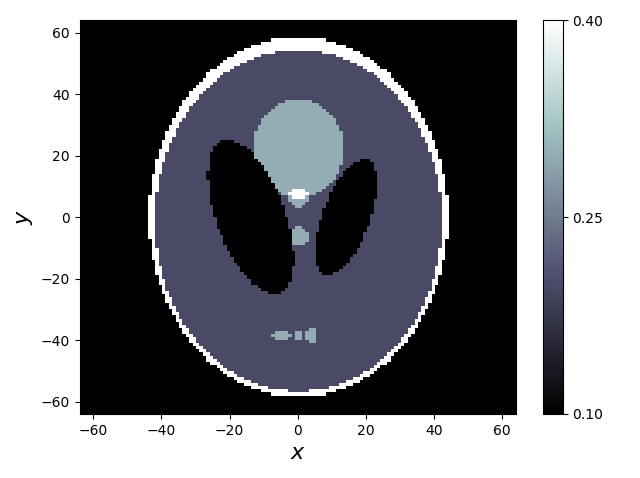}
		\caption{Phantom}
	\end{subfigure}
	\begin{subfigure}[t]{.49\linewidth}
		\includegraphics[width=\linewidth, trim={21mm 16mm 32mm 6mm}, clip]{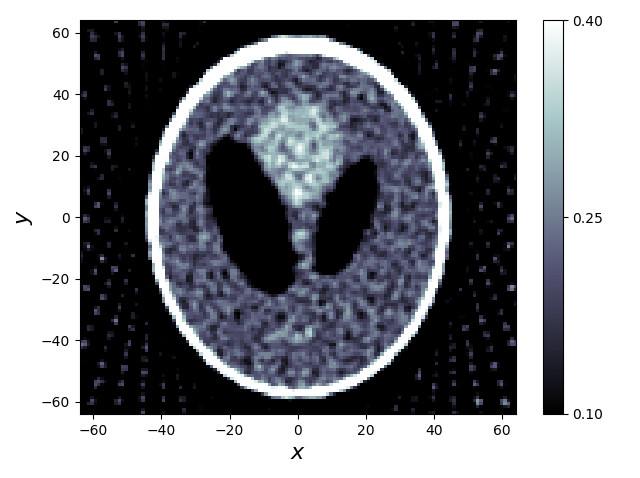}
		\caption{\ac{FBP}}
	\end{subfigure}
        \\[1em]	
	\begin{subfigure}[t]{.49\linewidth}
		\includegraphics[width=\linewidth, trim={21mm 16mm 32mm 6mm}, clip]{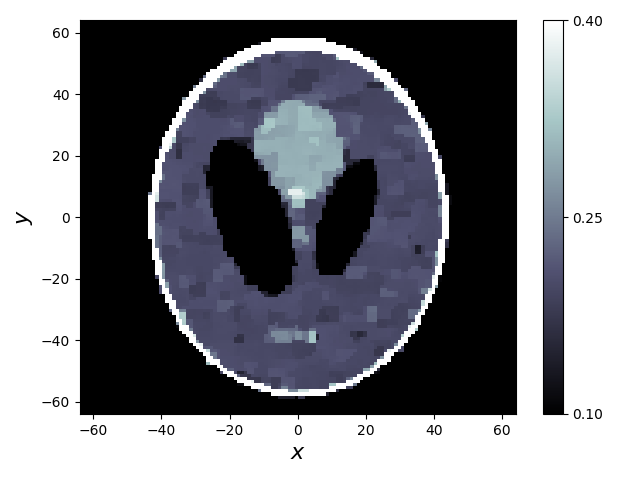}
		\caption{\ac{TV}}
	\end{subfigure}
	\begin{subfigure}[t]{.49\linewidth}
		\includegraphics[width=\linewidth, trim={21mm 16mm 32mm 6mm}, clip]{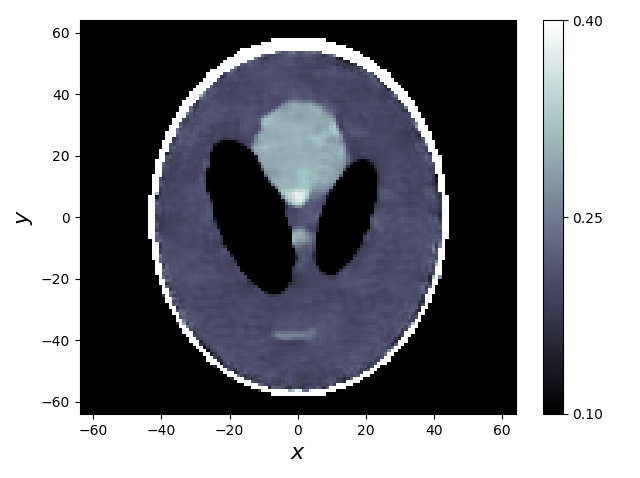}
		\caption{Partially learned gradient scheme}
	\end{subfigure}
	\caption{Reconstructing Shepp-Logan phantom using \ac{FBP}, \ac{TV} and the partially learned gradient scheme.
	  Data is simulated from the Shepp-Logan phantom, which attains values between $[0,1]$. All images 
	  are shown using a window set to $[0.1, 0.4]$ for improved contrast.}
	\label{fig:results_shepp}
\end{figure}

\begin{figure}[h]
	\centering	
	\begin{subfigure}[t]{.49\linewidth}
		\includegraphics[width=\linewidth, trim={23mm 16mm 27mm 6mm}, clip]{head_phantom}
		\caption{Phantom}
	\end{subfigure}
	\begin{subfigure}[t]{.49\linewidth}
		\includegraphics[width=\linewidth, trim={23mm 16mm 27mm 6mm}, clip]{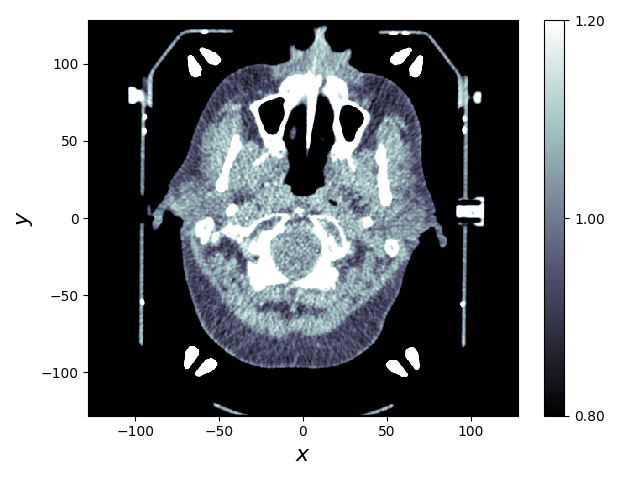}
		\caption{\ac{FBP}}
	\end{subfigure}
        \\[1em]	
	\begin{subfigure}[t]{.49\linewidth}
		\includegraphics[width=\linewidth, trim={23mm 16mm 27mm 6mm}, clip]{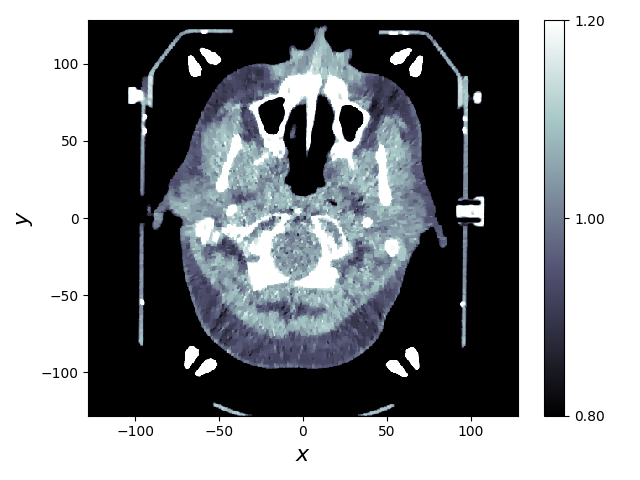}
		\caption{\ac{TV}}
	\end{subfigure}
	\begin{subfigure}[t]{.49\linewidth}
		\includegraphics[width=\linewidth, trim={23mm 16mm 27mm 6mm}, clip]{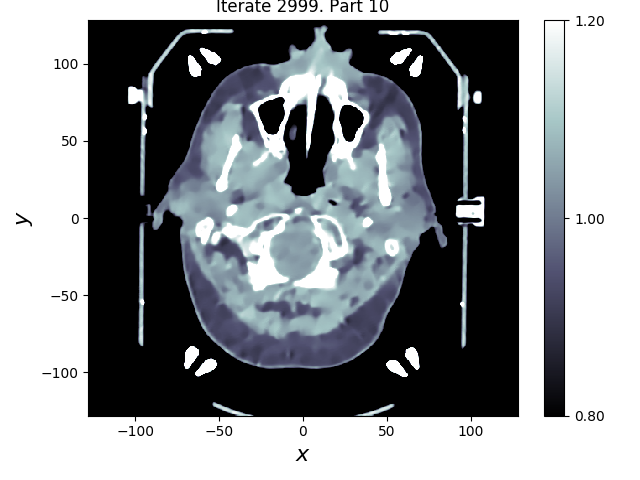}
		\caption{Partially learned gradient scheme}
	\end{subfigure}
	\caption{Reconstructing a head phantom using \ac{FBP}, \ac{TV} and the partially learned gradient scheme.
	  Data is simulated from a physiological head phantom, whuch includes some weak streaks (so these are part of the 
	  ground truth). All images are shown using a window set to $[-200, 200]$ \hounsfield.}
	\label{fig:results_head}
\end{figure}

\begin{figure}[h]
	\centering	
	\begin{subfigure}[t]{.32\linewidth}
		\includegraphics[width=\linewidth, trim={21mm 16mm 32mm 3mm}, clip]{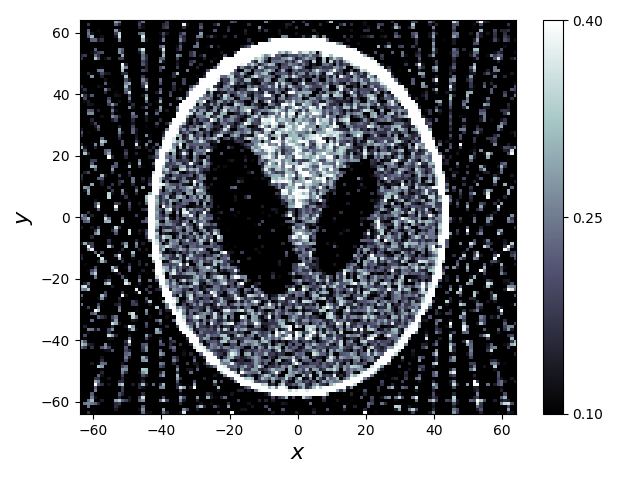}
		\caption{Initial}
		\label{fig:shepp_logan_result_iterate_0}
	\end{subfigure}
	\begin{subfigure}[t]{.32\linewidth}
		\includegraphics[width=\linewidth, trim={21mm 16mm 32mm 3mm}, clip]{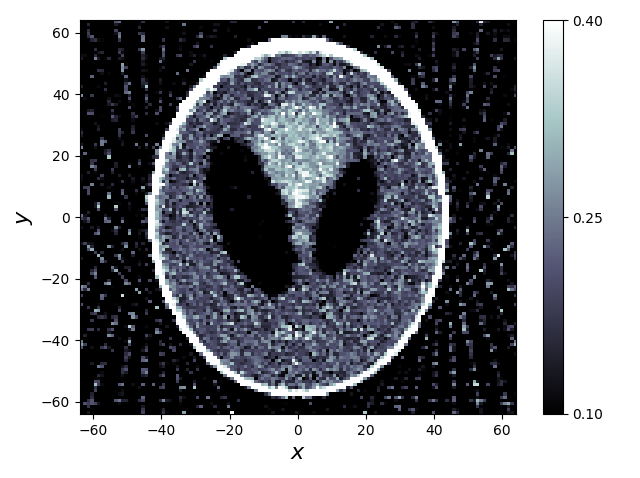}
		\caption{Iterate 1}
		\label{fig:shepp_logan_result_iterate_1}
	\end{subfigure}
	\begin{subfigure}[t]{.32\linewidth}
		\includegraphics[width=\linewidth, trim={21mm 16mm 32mm 3mm}, clip]{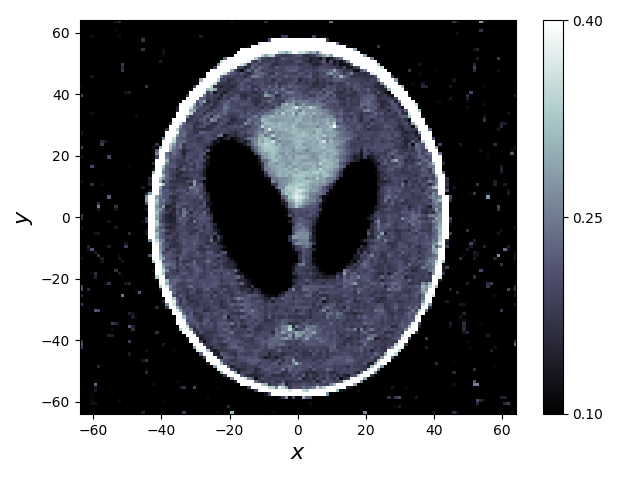}
		\caption{Iterate 2}
		\label{fig:shepp_logan_result_iterate_2}
	\end{subfigure}
        \\[1em]	
	\begin{subfigure}[t]{.32\linewidth}
		\includegraphics[width=\linewidth, trim={21mm 16mm 32mm 3mm}, clip]{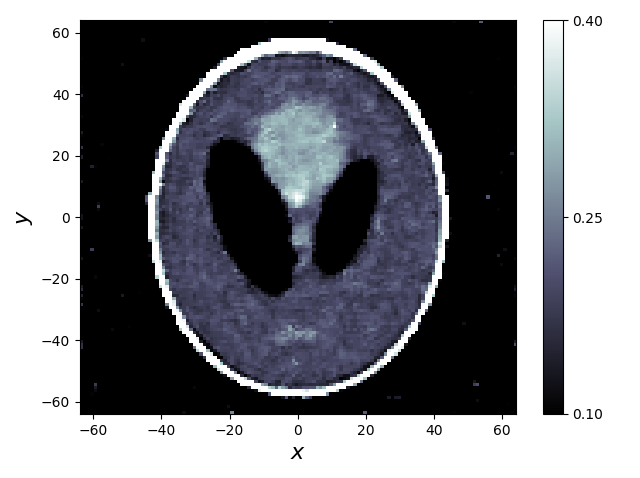}
		\caption{Iterate 3}
		\label{fig:shepp_logan_result_iterate_3}
	\end{subfigure}
	\begin{subfigure}[t]{.32\linewidth}
		\includegraphics[width=\linewidth, trim={21mm 16mm 32mm 3mm}, clip]{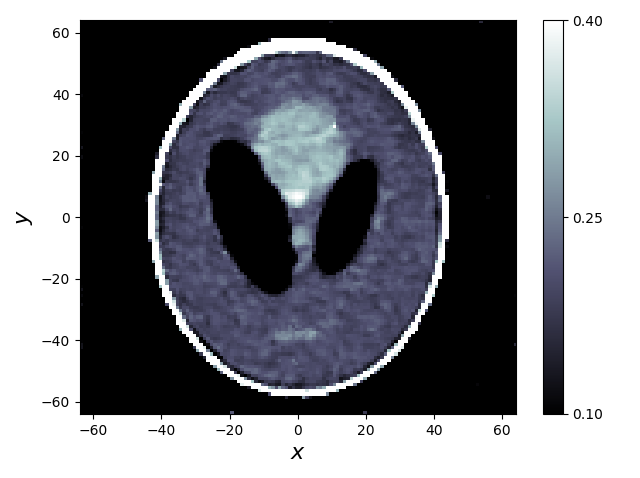}
		\caption{Iterate 4}
		\label{fig:shepp_logan_result_iterate_4}
	\end{subfigure}
	\begin{subfigure}[t]{.32\linewidth}
		\includegraphics[width=\linewidth, trim={21mm 16mm 32mm 3mm}, clip]{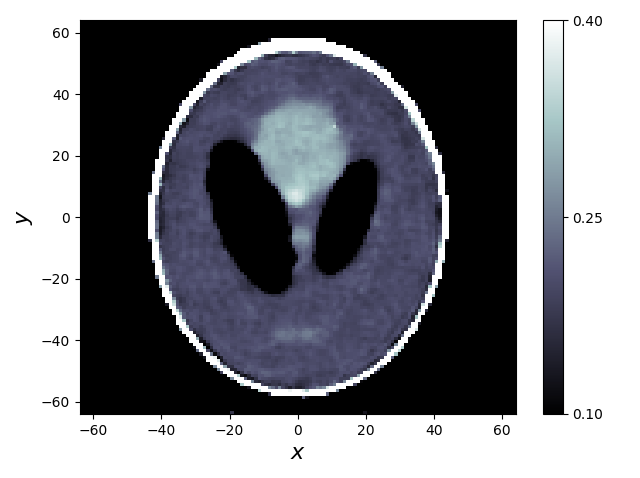}
		\caption{Iterate 5}
		\label{fig:shepp_logan_result_iterate_5}
	\end{subfigure}
        \\[1em]	
	\begin{subfigure}[t]{.32\linewidth}
		\includegraphics[width=\linewidth, trim={21mm 16mm 32mm 3mm}, clip]{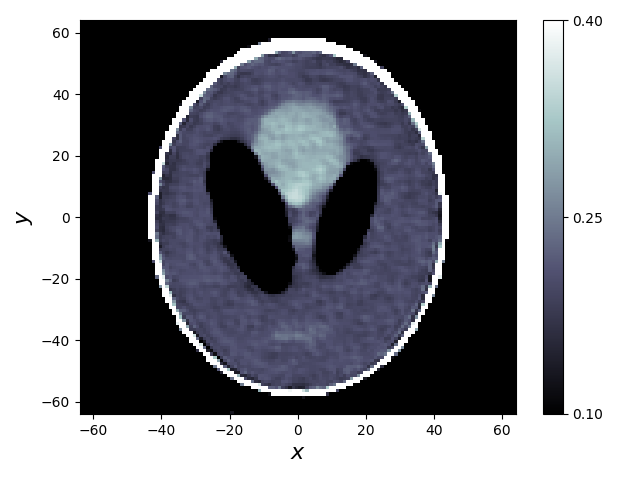}
		\caption{Iterate 6}
		\label{fig:shepp_logan_result_iterate_6}
	\end{subfigure}
	\begin{subfigure}[t]{.32\linewidth}
		\includegraphics[width=\linewidth, trim={21mm 16mm 32mm 3mm}, clip]{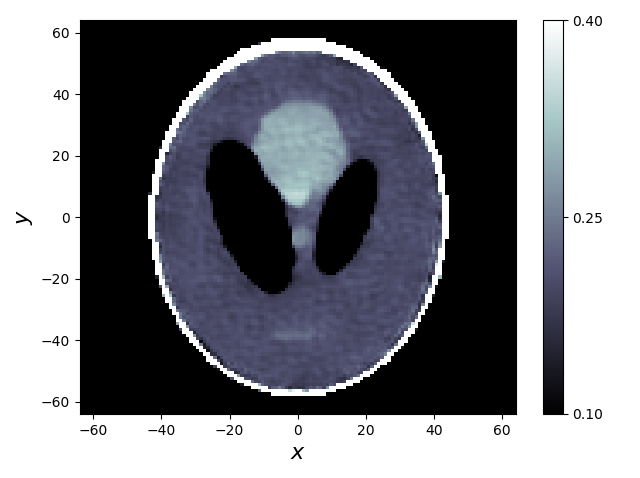}
		\caption{Iterate 8}
		\label{fig:shepp_logan_result_iterate_8}
	\end{subfigure}
	\begin{subfigure}[t]{.32\linewidth}
		\includegraphics[width=\linewidth, trim={21mm 16mm 32mm 3mm}, clip]{shepp_logan_result_windowed_iterate_10}
		\caption{Iterate 10}
		\label{fig:shepp_logan_result_iterate_10}
	\end{subfigure}
	\caption{Iterates of the partially learned gradient scheme when applied to reconstruct the Shepp-Logan phantom.
	  The initial iterate is given by the \ac{FBP} method.}
	\label{fig:iterates}
\end{figure}

We compare the reconstructions of the partially learned algorithm with the \ac{FBP} and \ac{TV} reconstructions for both the ellipse and head datasets and computed the \ac{PSNR}, runtime and performed a visual comparison. The quantitative results are given in \cref{tab:quant_res}, which visualizations are available in \cref{fig:results_shepp,fig:results_head}. We also display some partial results of the iterative algorithm for reference in \cref{fig:iterates}.

We note that for the ellipse data, the \ac{FBP} algorithm performs very poorly under the high noise while the \ac{TV} and learned methods give comparable results. This is expected given that \ac{TV} regularization is very competitive for images of this type, nonetheless the learned method does outperform the \ac{TV} method by approximately \unit{2}{\decibel} and the visual result looks slightly more appealing, with less randomly occurring structures and significantly less stair-casing. 

For the head dataset where the noise is lower, the filtered back-projection reconstruction performs much better, and is arguably comparable to the \ac{TV} regularized reconstruction. The learned reconstruction provides (perhaps too) smooth images, where we note that especially in the boundary regions, e.g. in the air-skin boundaries and around the bone the algorithm performs amiably.

In addition to the visual results, we see that the learned algorithm is significantly better than the \ac{TV} reconstruction w.r.t the \ac{PSNR}, giving an $>$\unit{5}{\decibel} improvement. Finally, the runtime of the algorithm, while being slightly slower than traditional filtered-back-projection, is significantly faster than the \ac{TV} method as shown in \cref{tab:quant_res}.

\begin{table}[ht]
\centering
\begin{tabular}{ l rr rr}
\toprule
 & \multicolumn{2}{c}{\ac{PSNR} (\decibel)} & \multicolumn{2}{c}{Runtime (\milli\second)} \\
\cmidrule(lr){2-3}
\cmidrule(lr){4-5}
Method       & Ellipses & Heads & Ellipses & Heads \\
\midrule
\ac{FBP}  & 19.75 & 36.12 & \textbf{4} & \textbf{130} \\
Learned & \textbf{32.02} & \textbf{43.82} & 58 & 430\\
\ac{TV} & 29.83 & 38.40 & 11\,963 & 173\,845\\
\bottomrule
\end{tabular}
\caption{Comparison of the learned method with standard methods.}
\label{tab:quant_res}
\end{table}

\subsection{Impact of including gradient mappings}
The impact of including the gradient mappings 
\[ 
\grad \big[ \LogLikelihood\bigl(\ForwardOp(\cdot), \data\bigr) \big], \grad\RegOp \colon \RecSpace \to \RecSpace
\]
in the partially learned gradient scheme can be empirically analysed. We do this by training the network in the exact same manner with and without the gradients added and then performing 100 reconstructions of the Shepp-Logan phantom with the respective methods. 

\begin{figure}
	\centering	
	\begin{subfigure}[t]{.49\linewidth}
		\includegraphics[width=\linewidth, trim={21mm 16mm 27mm 6mm}, clip]{shepp_logan_phantom_windowed}
		\caption{Phantom}
	\end{subfigure}
	\begin{subfigure}[t]{.49\linewidth}
		\includegraphics[width=\linewidth, trim={21mm 16mm 27mm 6mm}, clip]{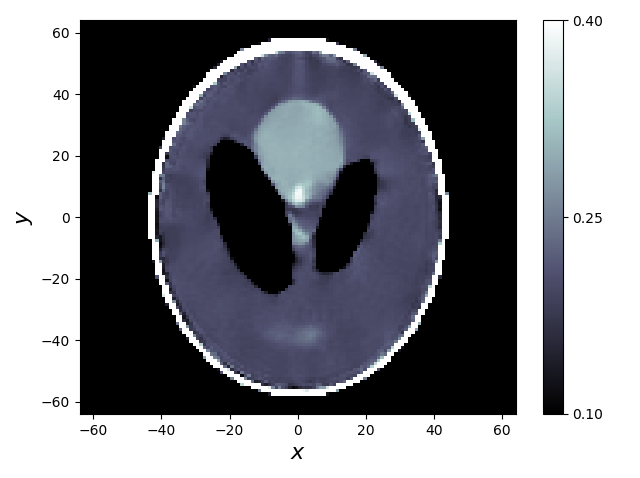}
		\caption{No gradient}
	\end{subfigure}
	\\[1em]
	\begin{subfigure}[t]{.49\linewidth}
		\includegraphics[width=\linewidth, trim={21mm 16mm 27mm 6mm}, clip]{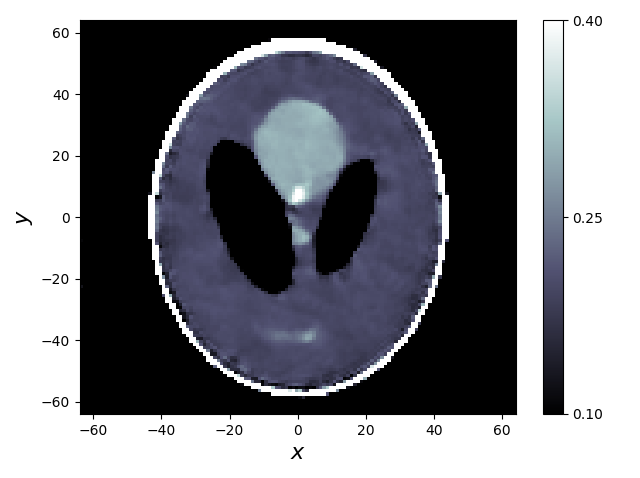}
		\caption{Only gradient of data discrepancy}
	\end{subfigure}
	\begin{subfigure}[t]{.49\linewidth}
		\includegraphics[width=\linewidth, trim={21mm 16mm 27mm 6mm}, clip]{shepp_logan_result_windowed_iterate_10}
		\caption{Gradient of both data discrepancy and regularizer}
	\end{subfigure}
	\caption{Comparison of reconstructions using the partially learned gradient scheme with and without the gradient 
	information. Note that gradient of data discrepancy includes the derivative of the forward operator.}
	\label{fig:results_gradients}
\end{figure}

Without the gradients, the \ac{PSNR} was \unit{29.65}{\decibel} while it was \unit{30.51}{\decibel} with the gradient of the data discrepancy. This can be compared to \unit{32.02}{\decibel} with both gradients. Visual inspection also indicates that adding the gradients provides a sharper reconstruction with more detail, where especially in the case of no gradients the small inserts are barely visible. Performance wise, the method took \unit{19}{\milli\second} without the gradients, \unit{64}{\milli\second} with the gradient of the data discrepancy and \unit{66}{\milli\second} with both gradients. See \cref{fig:results_gradients} for a visual comparison.

\section{Discussion}\label{sec:Disc}
The partially learned gradient scheme differs significantly from the current paradigm for regularization of inverse problems, so there are several remarks that deserve a closer discussion. 

\subsection{Theory}
The partially learned gradient scheme is presented with a strong emphasis on the algorithmic aspects, its implementation and its performance. There are however several interesting theoretical issues that deserve closer attention.

\paragraph{Deep neural networks in function spaces}
The scheme in \cref{alg:main_alg} is formulated in a functional analytic setting. However, the theory for deep neural networks is not well established in the infinite dimensional setting and there are several open issues that remain to be answered, such as determining what class of operators can be approximated by a given deep neural network and to what accuracy \cite{Gu16}. 

Another aspect relates to usage of probabilistic notions in infinite dimensional vector spaces. The classical theory for deep learning deals with finite, or at most countable data where the law of large numbers holds. In the infinite dimensional setting one needs to be more careful regarding which topologies that are used. It is clearly advantageous to work with spaces where one can prove various forms of weak convergence of probability measures and state and prove results corresponding to the law of large numbers, see \cite{VaTaCh87,St10}. Hence, applying deep learning to infinite dimensional spaces is associated with a number of fundamental questions regarding convergence of the learning, and if it converges, in what sense?

\paragraph{Regularizing properties}
A theoretical topic of interest is to prove that the given reconstruction scheme constitutes a formal regularization in the sense of \cite{ScGrGrHaLe09}, that is proving existence, stability and convergence. 

\begin{figure}
	\centering	
	\begin{subfigure}[t]{.49\linewidth}
		\includegraphics[width=\linewidth, trim={21mm 16mm 27mm 6mm}, clip]{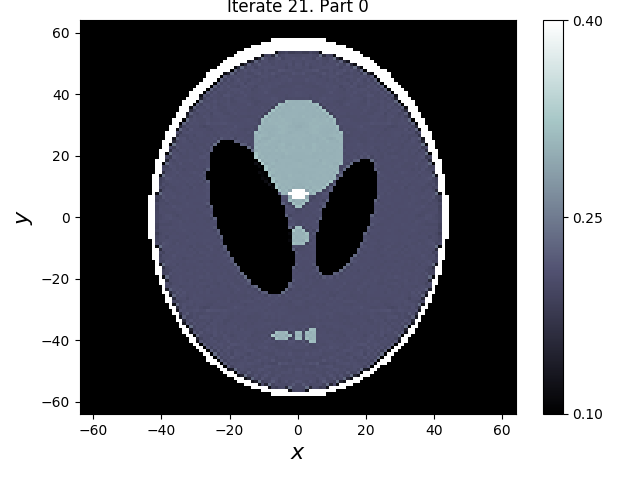}
		\caption{FBP}
	\end{subfigure}
	\begin{subfigure}[t]{.49\linewidth}
		\includegraphics[width=\linewidth, trim={21mm 16mm 27mm 6mm}, clip]{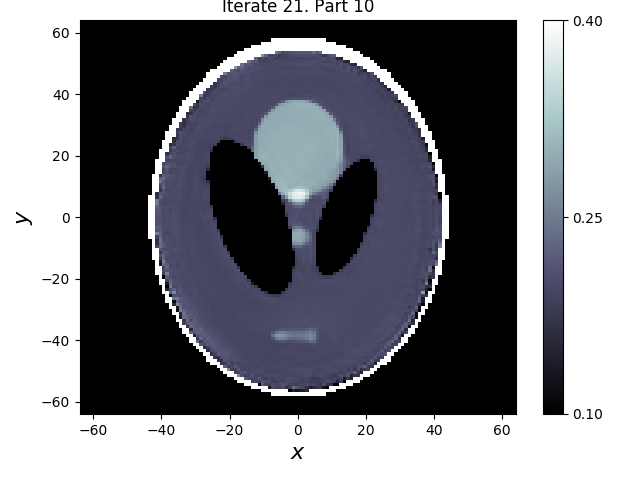}
		\caption{Learned}
	\end{subfigure}
	\caption{Example of partially learned reconstruction when applied to fully sampled and noiseless data while the parameter $\vparam$ is trained on sparsely sampled and noisy data.}
	\label{fig:learned_no_noise}
\end{figure}

Existence for $\ForwardOpInvLearned \colon \DataSpace \to \RecSpace$ is a non-issue since it this operator is given by a finite number of compositions of well-defined operators. Next, assume that the gradients of the data discrepancy 
$\LogLikelihood\bigl(\ForwardOp(\cdot), \data\bigr) 
\colon \RecSpace \to \Real$
and the regularizer $\RegOp \colon \RecSpace \to \Real$ are Lipschitz continuous. Then the partially learned pseudo-inverse $\ForwardOpInvLearned$ is also Lipschitz continuous, which in turn implies stability. The final consideration concerns convergence, which is formally defined as
\[
	\Bigl\lVert \ForwardOpInvLearned\bigl(\ForwardOp(\signaltrue) + \noise\bigr) - \signal^* \Bigr\rVert_\RecSpace \to 0
	\quad
	\text{whenever}
	\quad
	\lVert \noise \rVert_\DataSpace \to 0 
\]
for some parameter choice rule for the hyper-parameters in \cref{sec:hyperparam} and the training data, which uniquely define $\vparam$, and where $\signal^* \in \RecSpace$ is some minimum norm solution to \cref{eq:InvProb}. Clearly, the above convergence criteria can only be satisfied in general if the hyper-parameters and training data used for learning $\vparam$ are re-chosen as the data noise tends to zero. To conclude, in \cref{fig:learned_no_noise} we compare the result of applying the partially learned reconstruction scheme to noiseless data while training the parameter $\vparam$ against noisy data. It is reasonable to expect a significantly lower error if the method was re-trained on noiseless data, but we are currently unable to give a rigorous proof that this would converge to zero.

\subsection{Use cases}\label{sec:UseCases}
The framework for partially learned reconstruction was primarily motivated by a number of use cases involving ill-posed inverse problems. A number of challenges naturally arise when classical regularization is applied to solve the associated inverse problems and below we describe how these challenges can be resolved using a partially learned reconstruction scheme.
	
\paragraph{Computational feasibility}
The forward operator is an important part of a regularization and the more accurately it models the relation between signal and data, the better the outcome. Usage of more accurate forward models is however almost always computationally more demanding. Likewise, more elaborate regularization schemes that are better at utilizing the available a priori knowledge are often also computationally more demeaning. As an example, several of the more advanced regularizers in the literature exploit some kind of sparsity using a $L_1$-like norm \cite{ScGrGrHaLe09}. Such regularizers typically give rise to non-differentiable objective functional that require using optimization algorithms from non-smooth analysis for their efficient solution. Finally, there may also be reconstruction parameters that, unlike regularization parameter(s), do not influence the reconstruction quality. The role of the reconstruction parameters is to ensure the method is as efficient as possible, so these affect the speed of reconstruction.
	
Computational feasibility becomes especially critical in imaging applications since these involve very large-scale data structures. In such setting, variational and iterative regularization schemes quickly become infeasible even for applications with moderate time requirements despite usage of state of the art algorithms. 
	
The learned method \cref{alg:main_alg} improves upon this by having an a priori defined run-time which can be tweaked by using more or less iterates or a more complicated updating operator. By learning, we thus learn a optimal reconstruction scheme for a given execution time. Note that the run-time of the method on our examples is significantly faster than the \ac{TV} regularized method.
	
\paragraph{Nuisance parameters}
Nuisance parameters are additional unknowns that need to be reconstructed alongside the signal. They are not of primary interest, but they nevertheless need to be reconstructed. As an example, in certain tomographic applications the acquisition geometry (sampling of the data manifold $\datamanifold$) is partially unknown, so the nuisance parameters would be those needed for a precise description of said geometry. Another is use of a more accurate forward model, which often introduces nuisance parameters.
	
A common approach is to adopt an intertwined scheme in which each iterate involves updating the signal by reconstructing it from data using the  previous value for the nuisance parameter(s), followed by updating the nuisance parameter(s) by reconstructing them from data and making use of the previous value of recently updated signal. 

Our learned reconstruction scheme \cref{alg:main_alg} can easily be extended to include such intertwined schemes.
	
\paragraph{Regularization parameter selection rule}
Regularization parameter(s) govern the balancing between preventing over-fitting against the need to have a solution that generates data, which  is consistent with measurements. To have an appropriate parameter choice rule is critical for success.
	
Unfortunately, there is little theory to guide how to choose the regularization parameter(s). Mathematical results often study asymptotic behaviour of a parameter choice rule as data noise level tends to zero. Results mainly cover the case when noise in data is additive Gaussian and its magnitude can be reliably estimated \cite{EnHaNe00}, even though there are extensions for other noise types as well. Nevertheless, many of these assumptions are often not met in reality.

Some work has been done in selecting an optimal parameter using learning \cite{CaCaReScVa15}. The proposed method encompasses this since the regularization parameter (and other optimization related parameters) are included in the learned updating operator and thus optimally selected from the training data.
	
\paragraph{Feature reconstruction}
Reconstructing the signal is in many applications merely one part of a more elaborate scheme of transforming measured data to knowledge. As an example, in tomographic imaging the reconstructed image serves as input for an image analysis part. The latter often involves complex procedures, like segmentation and object recognition, that currently require involvement of human expertise.  

There is a growing trend in including some of these into the inverse problem that is referred to as feature reconstruction. To some extent, compressed sensing can be seen as an example of feature reconstruction where the sparse coding dictionary is the feature extraction part. Other examples are joint image reconstruction and segmentation \cite{RaRi07,BaSi11,Lo11} and shape based reconstruction \cite{FiGrSc12,GoXuReOkSu12,OkChDoRaBa17}. Such feature reconstruction methods are however hard to analyse theoretically and current methods are limited in the type of feature extraction capabilities they can include. They also tend to be computationally demanding.

It is natural to perform feature reconstruction by adding a feature extraction network to the learned reconstruction scheme such as in \cite{DiSiBoWeHe17}. The proposed framework could in a similar way be extended to feature reconstruction by composing the learned reconstruction operator $\ForwardOpInvLearned$ with a feature extraction operator $\FeatureExt \colon \RecSpace \to \FeatureSpace$ where $\FeatureSpace$ is a vector space of features. If the latter is differentiable, which is the case for deep learning based feature extractors, then we can define the loss functional \cref{eq:LossAbstract} using the composed operator $\FeatureExt \circ \ForwardOpInvLearned$. This allows for \emph{truly end-to-end optimization of task dependent reconstruction schemes for general inverse problems}.

\subsection{Stability}
A general question often asked when learning is applied to some problem is whether the method generalizes to other problems, e.g. if a method that is trained on a specific dataset can be applied to another dataset or to what extent one can change to forward operator without re-training.

Note first that the partially learned gradient scheme does not have an explicit regularization parameter, instead its regularization properties are implicitly contained in the training dataset (and to some extent in the hyper-parameters). Hence, a significant change in the training dataset (notably, a change of scaling) would require a re-training. On the other hand, empirical numerical experience suggests that dependence is relatively weak, at least for the tomographic reconstruction problems we considered. Specifically, we were able to successfully pre-train the system using a simplified acquisition geometry, a linearised forward operator, different domain size and significantly simplified phantoms and then successfully use this to train the network for the much more complicated heads dataset.

Finally, numerical experiments also suggests that changing the forward operator requires only a modest fine-tuning where the given parameters $\vparam$ can be used as an initial guess.

\section{Conclusion and future work}\label{sec:Conc}
We have presented a partially learned approach for solving ill-posed inverse problems that can integrate prior knowledge about the inverse problem with learning from training data. The presented method works with any non-linear operator and the method could easily be applied to a wide range of problems. Numerical experiments on tomographic data shows that the method gives notably better reconstructions than traditional \ac{FBP} and \ac{TV} regularization. Furthermore, adding prior information improves the reconstruction. In conclusion, using prior knowledge about the forward operator, data acquisition, data noise model and regularizer can significantly improve the performance of deep learning based approaches for solving inverse problems, and especially so when the available training data is much smaller than the size of the parameter space.

An obvious next step is to tackle fully three-dimensional tomographic problems while training on two-dimensional datasets. It would also be of interest to improve upon the choice of regularizer by adding more regularizers and/or more advanced regularizers such as wavelet based regularizers. Other more elaborate extensions are outlined below.

\paragraph{Extension to other iterative schemes}
The given iterative method is based of the gradient descent scheme, but this scheme is known to by sub-optimal in the case of non-differentiable objective functions. A natural extension of the scheme is thus to instead consider iterative schemes better suited for this use case. One such iterative scheme is the (non-linear) Chambolle-Pock algorithm \cite{ChPo10,Va14} for solving problems of the form
\[
\min_{\primal \in \RecSpace} \Bigl[ \OpF\bigl( \OpK(\primal)\bigr) + \OpG(\primal) \Bigr]
\]
where $\OpK \colon \RecSpace \to \ProdSpace$ is a (possibly non-linear) operator between Banach spaces $\RecSpace$ and $\ProdSpace$. The scheme is given by \cref{alg:cp} and the \emph{proximal} operators in \cref{alg:cp}�are given by
\begin{align*}
\ProxOp_{\sigma \OpF^*}(\dual) =& \argmin_{\dual' \in \ProdSpace} \Bigl[�\OpF^*(\dual') + \frac{1}{2 \sigma} \bigl\Vert \dual' - \dual \bigr\Vert_\ProdSpace^2 \Bigr] \\
\ProxOp_{\tau \OpG}(\primal) =& \argmin_{\primal' \in \RecSpace} \Bigl[ \OpG(\primal') + \frac{1}{2 \tau} \bigl\Vert \primal' - \primal \bigr\Vert_\RecSpace^2 \Bigr]
\end{align*}
where $\OpF^*$ is the Fenchel conjugate of $\OpF$. The special case of \ac{TV} regularized reconstruction for \eqref{eq:InvProb} amounts to selecting
\[
\OpK \colon \RecSpace \to \DataSpace \times \RecSpace^d
\quad
\text{as}
\quad
\OpK(\primal) := \bigl[ \ForwardOp(\primal), \nabla \primal \bigr]
\]
where $d$ is the dimension of the space and 
\[�\OpF\bigl([y_1, y_2] \bigr) := \lVert y_1 - \data \rVert_2^2 + \lVert y_2 \rVert 
   \quad\text{and}\quad
   \OpG(\primal) := 0.
\]
The resulting algorithm is summarized in \cref{alg:cp}.
\begin{algorithm}
	\caption{Non-linear Chambolle-Pock algorithm}\label{alg:cp}
	\begin{algorithmic}[1]
		\State Given: $\sigma, \tau > 0$ s.t. $\sigma \tau \lVert \OpK \rVert^2 < 1$, $\theta \in [0, 1]$ and $\primal_0 \in \RecSpace$, $\dual_0 \in \ProdSpace$. 
		\For{$i = 1, \dots, I$}
		\State $\dual^{i + 1} = \ProxOp_{\sigma \OpF^*}\bigl(\dual^{i} + \sigma \OpK(\bar{\primal}^{i}) \bigr)$ 
		\State $\primal^{i + 1} = \ProxOp_{\tau \OpG}\bigl(\primal^{i} - \tau [ \partial \OpK(\primal^{i})]^* (\dual^{i + 1}) \bigr)$
		\State $\bar{\primal}^{i + 1} = \primal^{i + 1} + \theta (\primal^{i + 1} - \primal^{i})$
		\EndFor
	\end{algorithmic}
\end{algorithm}

To introduce a learning component, one may either learn the primal proximal ($\ProxOp_{\tau \OpG}$) or the dual proximal $(\ProxOp_{\sigma \OpF^*})$, or both. Some recent papers have approached learning the primal proximal operator \cite{ChLiPoViSa17, YaSuLiXu16} in the scope of \ac{ADMM}, but these do not consider learning the dual. It is likely that learning the dual proximal offers an advantage since this allows the inclusion of various operators into the learning. To illustrate this, one can learn a proximal operator for the directional wavelet coefficients of a signal. This is successfully done for de-noising \cite{WaveNet}, and would likely by useful for reconstruction as well.

\begin{figure}
	\centering	
	\begin{subfigure}[t]{.49\linewidth}
		\includegraphics[width=\linewidth, trim={21mm 16mm 27mm 6mm}, clip]{shepp_logan_phantom_windowed}
		\caption{Phantom}
	\end{subfigure}
	\begin{subfigure}[t]{.49\linewidth}
		\includegraphics[width=\linewidth, trim={21mm 16mm 27mm 6mm}, clip]{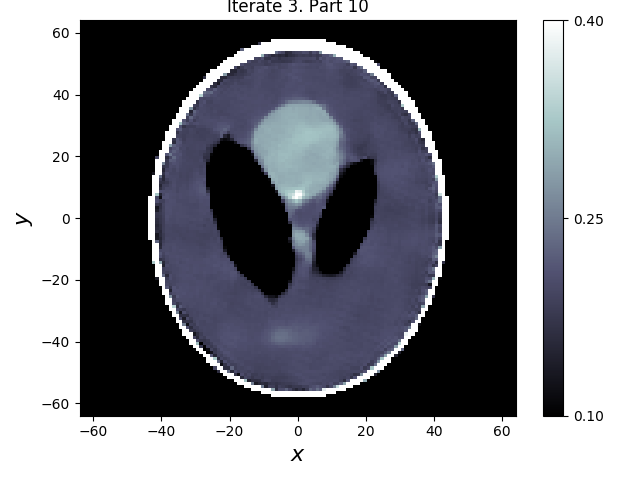}
		\caption{Reconstruction with learned proximal}
	\end{subfigure}
	\caption{Reconstruction with a partially learned dual proximal using the method in \cref{alg:cp}.}
	\label{fig:learned_dual}
\end{figure}

Learning the dual proximal also allows one to incorporate memory into the algorithm. This can be done for the above case by extending the operator $\OpK$ so that it also contains a zero component:
\[
\OpK \colon \RecSpace \to \DataSpace \times \RecSpace^d \times \RecSpace^M
\quad\text{with}\quad
\OpK(\primal) := \bigl[ \ForwardOp(\primal), \nabla \primal, 0 \bigr].
\]
We intend to further elaborate on this approach in an upcoming paper, at this stage we settle with providing an example reconstruction shown in \cref{fig:learned_dual}.

\paragraph{Choice of discretization}
The given examples were performed using the simplest discretization of the space $\RecSpace$, a pixel basis, but  \cref{alg:main_alg} also works with other representations such as Fourier, wavelet or shearlet coefficients. This could in many cases be better suited for the inverse problem in question and especially so if the operator $\ForwardOp$ or the regularizer $\RegOp$ has a simple form in this representation, as with the Fourier transform for \ac{MRI}.

\paragraph{Choice of error functional}
We simply investigated the squared norm error function, $\errorfunc(\signal) = \lVert \signal - \signaltrue \rVert_X^2$, but experience tells us that this is perhaps not the best predictor of human observer performance on a given image. For example, the given algorithm gives a \unit{5}{\decibel} improvement over the \ac{TV} algorithm, but by visual inspection the improvement is not equally drastic. A possible way to improve this and further leverage the power of the learning approach is to use a more sophisticated error functional. Here, performing end-to-end optimization should be a feasible alternative, i.e. instead maximize some type of task based measure.

\section{Acknowledgements}
The work of both authors was supported by the Swedish Foundation of Strategic Research grant AM13-0049 and Industrial PhD grant ID14-0055. Jonas Adler's work was also supported by Elekta.

\bibliographystyle{plain}
\bibliography{learning_refs}

\begin{thebibliography}{10}

\bibitem{AbAgBaBrCh15}
M.~Abadi, A.~Agarwal, P.~Barham, E.~Brevdo, Z.~Chen, C.~Citro, G.~S. Corrado,
  A.~Davis, J.~Dean, M.~Devin, S.~Ghemawat, I.~Goodfellow, A.~Harp, G.~Irving,
  M.~Isard, J.~Yangqing, R.~Jozefowicz, L.~Kaiser, M.~Kudlur, J.~Levenberg,
  D.~Man{\'e}, R.~Monga, S.~Moore, D.~Murray, C.~Olah, M.~Schuster, J.~Shlens,
  B.~Steiner, I.~Sutskever, K.~Talwar, P.~Tucker, V.~Vanhoucke, V.~Vasudevan,
  F.~Vi{\'e}gas, O.~Vinyals, P.~Warden, M.~Wattenberg, M.~Wicke, Y.~Yu, and
  X.~Zheng.
\newblock {TensorFlow}: {L}arge-scale machine learning on heterogeneous
  systems.
\newblock ArXiv:cs.DC 1603.04467, ArXiv, 2015.

\bibitem{AdKoOk17}
J.~Adler, H.~Kohr, and O.~{\"O}ktem.
\newblock Operator discretization library {(ODL)}.
\newblock Software available from
  \href{https://github.com/odlgroup/odl}{https://github.com/odlgroup/odl},
  2017.

\bibitem{AnDeCoHoPf16}
M.~Andrychowicz, M.~Denil, S.~G. Colmenarejo, M.~W. Hoffman, D.~Pfau,
  T.~Schaul, and N.~de~Freitas.
\newblock Learning to learn by gradient descent by gradient descent.
\newblock ArXiv:cs.NE 1606.04474, ArXiv, 2016.
\newblock Report available from
  \href{https://arxiv.org/abs/1606.04474}{https://arxiv.org/abs/1606.04474}.

\bibitem{ArMaTsSt12}
M.~Argyrou, D.~Maintas, C.~Tsoumpas, and E.~Stiliaris.
\newblock Tomographic image reconstruction based on artificial neural network
  {(ANN)} techniques.
\newblock In {\em Nuclear Science Symposium and Medical Imaging Conference
  (NSS/MIC), 2012 IEEE}, 2012.

\bibitem{BaSi11}
K.~J. Batenburg and J.~Sijbers.
\newblock {DART}: {A} practical reconstruction algorithm for discrete
  tomography.
\newblock {\em IEEE Transactions on Image Processing}, 20(9):2542--2553, 2011.

\bibitem{BeLaZa08}
M.~Bertero, H.~Lant{\'e}ri, and L.~Zanni.
\newblock Iterative image reconstruction: a point of view.
\newblock In Y.~Censor, M.~Jiang, and A.~K. Louis, editors, {\em Proceedings of
  the Interdisciplinary Workshop on Mathematical Methods in Biomedical Imaging
  and Intensity-Modulated Radiation (IMRT), Pisa, Italy}, pages 37--63, 2008.

\bibitem{CaCaReScVa15}
L.~Calatroni, C.~Cao, J.~C. De~Los~Reyes, C.-B. Sch{\"o}nlieb, and T.~Valkonen.
\newblock Bilevel approaches for learning of variational imaging models.
\newblock ArXiv:math.OC 1505.02120, ArXiv, 2015.
\newblock To appear in RICAM special issue, report available from
  \href{https://arxiv.org/abs/1505.02120}{https://arxiv.org/abs/1505.02120}.

\bibitem{CaReSc16}
C.~Cao, J.~C. De~Los~Reyes, and C.-B. Sch{\"o}nlieb.
\newblock Learning optimal spatially-dependent regularization parameters in
  total variation image restoration.
\newblock ArXiv:math.OC 1603.09155, ArXiv, 2016.
\newblock To appear in Inverse Problems, report available from
  \href{https://arxiv.org/abs/1603.09155}{https://arxiv.org/abs/1603.09155}.

\bibitem{ChPo10}
A.~Chambolle and T.~Pock.
\newblock A first-order primal-dual algorithm for convex problems with
  applications to imaging.
\newblock Technical Report 00490826, HAL-archives, 2010.
\newblock Report available from
  \href{https://hal.archives-ouvertes.fr/hal-00490826}{https://hal.archives-ouvertes.fr/hal-00490826}.

\bibitem{ChLiPoViSa17}
J.~H.~R. Chang, C.-L. Li, B.~Poczos, B.~V.~K. Vijaya~Kumar, and A.~C.
  Sankaranarayanan.
\newblock One network to solve them all --- solving linear inverse problems
  using deep projection models.
\newblock ArXiv:cs.CV 1703.09912, ArXiv, 2017.
\newblock Report available from
  \href{https://arxiv.org/abs/1703.09912}{https://arxiv.org/abs/1703.09912}.

\bibitem{ReSc13}
J.~C. De~Los~Reyes and C.-B. Sch{\"o}nlieb.
\newblock Image denoising: Learning noise distribution via {PDE}-constrained
  optimisation.
\newblock {\em Inverse Problems and Imaging}, 7:1183--1214, 2013.

\bibitem{ReScVa16}
J.~C. De~Los~Reyes, C.-B. Sch{\"o}nlieb, and T.~Valkonen.
\newblock The structure of optimal parameters for image restoration problems.
\newblock {\em Journal of Mathematical Analysis and Applications},
  434:464--500, 2016.

\bibitem{ReScVa17}
J.~C. De~Los~Reyes, C.-B. Sch{\"o}nlieb, and T.~Valkonen.
\newblock Bilevel parameter learning for higher-order total variation
  regularisation models.
\newblock {\em Journal of Mathematical Imaging and Vision}, 57(1):1--25, 2017.

\bibitem{DiSiBoWeHe17}
S.~Diamond, V.~Sitzmann, S.~Boyd, G.~Wetzstein, and F.~Heide.
\newblock {Dirty Pixels:} {O}ptimizing image classification architectures for
  raw sensor data.
\newblock ArXiv:cs.CV 1701.06487, ArXiv, 2017.
\newblock Report available from
  \href{https://arxiv.org/abs/1701.06487}{https://arxiv.org/abs/1701.06487}.

\bibitem{EnHaNe00}
H.~W. Engl, M.~Hanke, and A.~Neubauer.
\newblock {\em Regularization of inverse problems}.
\newblock Number 375 in Mathematics and its Applications. Kluwer Academic
  Publishers, 2000.

\bibitem{FiGrSc12}
T.~Fidler, M.~Grasmair, and O.~Scherzer.
\newblock Shape reconstruction with a priori knowledge based on integral
  invariants.
\newblock {\em SIAM Journal of Imaging Sciences}, 5(2):726--745, 2012.

\bibitem{GoXuReOkSu12}
A.~Gopinath, G.~Xu, D.~Ress, O.~{\"O}ktem, S.~Subramaniam, and C.~Bajaj.
\newblock Shape-based regularization of electron tomographic reconstruction.
\newblock {\em IEEE Transactions on Medical Imaging}, 31(12):2241--2252, 2012.

\bibitem{Gu16}
W.~H. Guss.
\newblock Deep function machines: {G}eneralized neural networks for topological
  layer expression.
\newblock ArXiv:stat.ML 1612.04799, ArXiv, 2016.
\newblock Report available from
  \href{https://arxiv.org/abs/1612.04799}{https://arxiv.org/abs/1612.04799}.

\bibitem{HaWuPoMa17}
K.~Hammernik, T.~W{\"u}rfl, T.~Pock, and A.~Maier.
\newblock A deep learning architecture for limited-angle computed tomography
  reconstruction.
\newblock In K.~H. Maier-Hein, T.~M. Deserno, H.~Handels, and T.~Tolxdorff,
  editors, {\em Bildverarbeitung f{\"u}r die Medizin 2017: Algorithmen -
  Systeme - Anwendungen. Proceedings des Workshops vom 12. bis 14. M{\"a}rz
  2017 in Heidelberg}, pages 92--97. Springer-Verlag, Berlin, Heidelberg, 2017.

\bibitem{Ha97}
P.-C. Hansen.
\newblock {\em Rank-Deficient and Discrete Ill-Posed Problems: Numerical
  Aspects of Linear Inversion}, volume~4 of {\em SIAM Monographs on
  Mathematical Modeling and Computation}.
\newblock SIAM, 1997.

\bibitem{JiMcFrUn16}
K.~H. Jin, M.~T. McCann, E.~Froustey, and M.~Unser.
\newblock Deep convolutional neural network for inverse problems in imaging.
\newblock ArXiv:cs.CV 1611.03679, ArXiv, 2016.
\newblock Report available from
  \href{https://arxiv.org/abs/1611.03679}{https://arxiv.org/abs/1611.03679}.

\bibitem{KaNeSc08}
B.~Kaltenbacher, A.~Neubauer, and O.~Scherzer.
\newblock {\em Iterative Regularization Methods for Nonlinear Ill-posed
  Problems}, volume~6 of {\em {Radon Series on Computational and Applied
  Mathematics}}.
\newblock Walter de Gruyter, 2008.

\bibitem{WaveNet}
Eunhee Kang, Junhong Min, and Jong~Chul Ye.
\newblock Wavenet: a deep convolutional neural network using directional
  wavelets for low-dose x-ray {CT} reconstruction.
\newblock {\em CoRR}, abs/1610.09736, 2016.

\bibitem{LBFGS}
D.~C. Liu and J.~Nocedal.
\newblock On the limited memory bfgs method for large scale optimization.
\newblock {\em Math. Program.}, 45(3):503--528, December 1989.

\bibitem{Lo11}
A.~K. Louis.
\newblock Feature reconstruction in inverse problems.
\newblock {\em Inverse Problems}, 27:065010 (21pp), 2011.

\bibitem{relu}
Vinod Nair and Geoffrey~E. Hinton.
\newblock Rectified linear units improve restricted boltzmann machines.
\newblock In Johannes F{\"u}rnkranz and Thorsten Joachims, editors, {\em
  Proceedings of the 27th International Conference on Machine Learning
  (ICML-10)}, pages 807--814. Omnipress, 2010.

\bibitem{OkChDoRaBa17}
O.~{\"O}ktem, C.~Chen, N.~O. Domani{\c c}, P.~Ravikumar, and C.~Bajaj.
\newblock Shape-based image reconstruction using linearized deformations.
\newblock {\em Inverse Problems}, 33(3):035004 (33pp), 2017.

\bibitem{PaGiKaLoMa04}
P.~Paschalis, N.~D. Giokaris, A.~Karabarbounis, G.~K. Loudos, D.~Maintas, C.~N.
  Papanicolas, V.~Spanoudaki, Ch. Tsoumpas, and E.~Stiliaris.
\newblock Tomographic image reconstruction using artificial neural networks.
\newblock {\em Nuclear Instruments and Methods in Physics Research Section A:
  Accelerators, Spectrometers, Detectors and Associated Equipment},
  527(1--2):211--215, 2004.

\bibitem{PeBa13}
D.~M. Pelt and K.~J. Batenburg.
\newblock Fast tomographic reconstruction from limited data using artificial
  neural networks.
\newblock {\em IEEE Transactions on Image Processing}, 22(12):5238--5251, 2013.

\bibitem{PuWe17}
P.~Putzky and M.~Welling.
\newblock Recurrent inference machines for solving inverse problems.
\newblock Submitted to ICLR 2017, Toulon, France, April 24--26, 2017. Report
  available from
  \href{https://openreview.net/pdf?id=HkSOlP9lg}{https://openreview.net/pdf?id=HkSOlP9lg},
  2017.

\bibitem{RaRi07}
Ronny Ramlau and Wolfgang Ring.
\newblock A {Mumford--Shah} level-set approach for the inversion and
  segmentation of x-ray tomography data.
\newblock {\em Journal of Computational Physics}, 221(2):539 -- 557, 2007.

\bibitem{ScGrGrHaLe09}
O.~Scherzer, M.~Grasmair, H.~Grossauer, M.~Haltmeier, and F.~Lenzen.
\newblock {\em {Variational Methods in Imaging}}, volume 167 of {\em {Applied
  Mathematical Sciences}}.
\newblock Springer-Verlag, New York, 2009.

\bibitem{Sc07}
T.~Schuster.
\newblock {\em {The Method of Approximate Inverse: Theory and Applications}},
  volume 1906 of {\em {Lecture Notes in Mathematics}}.
\newblock Springer Verlag, Heidelberg, 2007.

\bibitem{St10}
A.~M. Stuart.
\newblock Inverse problems: A {B}ayesian perspective.
\newblock {\em Acta Numerica}, pages 451--559, 2010.

\bibitem{VaTaCh87}
N.~N. Vakhania, V.~I. Tarieladze, and S.~A. Chobanyan.
\newblock {\em Probability Distributions on {B}anach Spaces}.
\newblock Mathematics and Its Applications (Soviet Series). Kluwer Academic
  Publishers, 1987.

\bibitem{Va14}
T.~Valkonen.
\newblock A primal-dual hybrid gradient method for nonlinear operators with
  applications to {MRI}.
\newblock {\em Inverse Problems}, 30(5):055012, 2014.

\bibitem{AaPaCaJaBl16}
W.~van Aarle, W.~J. Palenstijn, J.~Cant, E.~Janssens, Folkert Bleichrodt,
  A.~Dabravolski, J.~Beenhouwer, K.~J. Batenburg, and J.~Sijbers.
\newblock Fast and flexible {X}-ray tomography using the {ASTRA} toolbox.
\newblock {\em Optics Express}, 24(22):25129--25147, 2016.

\bibitem{WuGhChMa16}
T.~W{\"u}rfl, F.~C. Ghesu, V.~Christlein, and A.~Maier.
\newblock Deep learning computed tomography.
\newblock In S.~Ourselin, L.~Joskowicz, M.~Sabuncu, G.~Unal, and W.~Wells,
  editors, {\em MICCAI 2016: Medical Image Computing and Computer-Assisted
  Intervention -- MICCAI 2016}, volume 9902 of {\em Lecture Notes in Computer
  Science}, pages 432--440. Springer-Verlag, 2016.

\bibitem{YaSuLiXu16}
Y.~Yang, J.~Sun, H.~Li, and Z.~Xu.
\newblock Deep {ADMM-Net} for compressive sensing {MRI}.
\newblock In D.~D. Lee, M.~Sugiyama, U.~V. Luxburg, I.~Guyon, and R.~Garnett,
  editors, {\em Advances in Neural Information Processing Systems}, volume~29,
  pages 10--18. Curran Associates, 2016.
\newblock Report available from
  \href{http://papers.nips.cc/paper/6406-deep-admm-net-for-compressive-sensing-mri.pdf}{http://papers.nips.cc/paper/6406-deep-admm-net-for-compressive-sensing-mri.pdf}.

\end{thebibliography}

\end{document}